\documentclass[reqno,twoside,11pt]{amsart}
\usepackage{amsmath}
\usepackage{amsfonts}
\usepackage{amssymb}
\usepackage{verbatim}
\usepackage{epsfig}

\IfFileExists{epsf.def}{\input epsf.def}{\usepackage{epsf}}

\usepackage{mathpazo}\usepackage{mathrsfs}
\DeclareFontFamily{OT1}{eusb}{} \DeclareFontShape{OT1}{eusb}{m}{n} {<5> <6> <7> <8> <9> <10> <11> <12> <14.4> eusb10}{}
\DeclareMathAlphabet{\eusb}{OT1}{eusb}{m}{n}

\DeclareFontFamily{OT1}{eusm}{} \DeclareFontShape{OT1}{eusm}{m}{n} {<5> <6> <7> <8> <9> <10> <11> <12> <14.4> eusm10}{}
\DeclareMathAlphabet{\eusm}{OT1}{eusm}{m}{n}

\DeclareFontFamily{OT1}{eufm}{} \DeclareFontShape{OT1}{eufm}{m}{n} {<5> <6> <7> <8> <9> <10> <11> <12> <14.4> eufm10}{}
\DeclareMathAlphabet{\mathfrak}{OT1}{eufm}{m}{n}

\DeclareFontFamily{OT1}{fraktura}{}
\DeclareFontShape{OT1}{fraktura}{m}{n} {<5> <6> <7> <8> <9> <10> <11> <12> <13> <14.4> [1.1] eufm10}{}
\DeclareMathAlphabet{\fraktura}{OT1}{fraktura}{m}{n}

\DeclareFontFamily{OT1}{cmfi}{} \DeclareFontShape{OT1}{cmfi}{m}{n} {<5> <6> <7> <8> <9> <10> <11> <12> <13> <14.4> [0.9] cmfi10}{}
\DeclareMathAlphabet{\cmfi}{OT1}{cmfi}{b}{n}

\DeclareFontFamily{OT1}{cmss}{} \DeclareFontShape{OT1}{cmss}{m}{n} {<5> <6> <7> <8> <9> <10> <11> <12> <13> <14.4> cmss10}{}
\DeclareMathAlphabet{\cmss}{OT1}{cmss}{m}{n}
\newtheoremstyle{thm}{1.8ex}{1.8ex}{\itshape\rmfamily}{} {\bfseries\rmfamily}{}{2ex}{}

\newtheoremstyle{def}{1.8ex}{1.8ex}{\slshape\rmfamily}{} {\bfseries\rmfamily}{}{2ex}{}

\newtheoremstyle{rem}{1.8ex}{1.8ex}{\rmfamily}{} {\bfseries\rmfamily}{}{2ex}{}
\newenvironment{proofsect}[1] {\vspace{0.2cm}\noindent{\rmfamily\itshape#1.}}{\qed\vspace{0.15cm}}%{\newline\vspace{0.15cm}}

\theoremstyle{thm}
\newtheorem{theorem}{Theorem}[section]
\newtheorem{lemma}[theorem]{Lemma}
\newtheorem{proposition}[theorem]{Proposition}

\newtheorem*{Main Theorem}{Main Theorem.}

\newtheorem*{special theorem}{Lindeberg-Feller Theorem for Martingales}

\theoremstyle{def}

\theoremstyle{rem}
\newtheorem{remark}[theorem]{Remark}

\numberwithin{equation}{section}

\renewcommand{\section}{\secdef\sct\sect}
\newcommand{\sct}[2][default]{%
\refstepcounter{section}
\addcontentsline{toc}{section}{{\tocsection {}{\thesection}{\!\!\!\!#1\dotfill}}{}}
\vspace{0.7cm}
\centerline{\scshape\thesection.\ #1} \nopagebreak \vspace{0.2cm}}
\newcommand{\sect}[1]{%
\vspace{0.4cm} \centerline{\large\scshape\rmfamily #1}
\vspace{0.2cm}}

\renewcommand{\subsection}{\secdef\subsct\sbsect}
\newcommand{\subsct}[2][default]{\refstepcounter{subsection}
\addcontentsline{toc}{subsection}
{{\tocsection{\!\!}{\hspace{1.2em}\thesubsection}{\!\!\!\!#1\dotfill}}{}}
\nopagebreak\vspace{0.45\baselineskip} {\flushleft\bf
\thesubsection~\bf #1.~}
\\*[3mm]\noindent
\nopagebreak}
\newcommand{\sbsect}[1]{\vspace{0.1cm}\noindent
\textbf{#1.~}\vspace{0.1cm}}

\renewcommand{\subsubsection}{%
\secdef \subsubsect\sbsbsect}
\newcommand{\subsubsect}[2][default]{%
\refstepcounter{subsubsection} 
\addcontentsline{toc}{subsubsection}{{\tocsection{\!\!}
{\hspace{3.05em}\thesubsubsection}{\!\!\!\!#1\dotfill}}{}}
\nopagebreak
\vspace{0.15\baselineskip} \nopagebreak {\flushleft\rmfamily
\itshape\thesubsubsection
\ \rmfamily #1\/.}\ }
\newcommand{\sbsbsect}[1]{\vspace{0.1cm}\noindent
\rmfamily \itshape
\arabic{section}.\arabic{subsection}.\arabic{subsubsection} \
\sffamily #1\/.\ }

\renewcommand{\caption}[1]{%
\vglue0.5cm
\refstepcounter{figure}
\begin{minipage}{0.9\textwidth}\small {\sc Figure~\thefigure. }#1\end{minipage}}

\newcommand{\diam}{\operatorname{diam}}

\newcommand{\textd}{\text{\rm d}\mkern0.5mu}

\newcommand{\texte}{\text{\rm e}}

\newcommand{\1}{\operatorname{\sf 1}}

\newcommand{\HH}{\mathcal H}

\newcommand{\LL}{\mathcal L}

\newcommand{\B}{\mathbb B}

\newcommand{\E}{\mathbb E}

\newcommand{\BbbP}{\mathbb P}

\newcommand{\R}{\mathbb R}

\newcommand{\Z}{\mathbb Z}

\newcommand{\scrC}{\mathscr{C}}

\newcommand{\scrF}{\mathscr{F}}

\newcommand{\scrO}{\mathscr{O}}

\newcommand{\scrW}{\mathscr{W} }

\newcommand{\twoeqref}[2]{(\ref{#1}--\ref{#2})}
\newcommand{\cc}{{\text{\rm c}}}

\def\myffrac#1#2 in #3{\raise 2.6pt\hbox{$#3 #1$}\mkern-1.5mu\raise 0.8pt\hbox{$#3/$}\mkern-1.1mu\lower 1.5pt\hbox{$#3 #2$}}
\newcommand{\ffrac}[2]{\mathchoice%
	{\myffrac{#1}{#2} in \scriptstyle}
	{\myffrac{#1}{#2} in \scriptstyle}
	{\myffrac{#1}{#2} in \scriptscriptstyle}
	{\myffrac{#1}{#2} in \scriptscriptstyle}
}

\newcommand{\DIV}{\,\text{\rm div}\,}
\newcommand{\GRAD}{\,\text{\rm grad}\,}

\newcommand{\hate}{\hat{\text{\rm e}}}
\newcommand{\hatv}{\hat{\text{\rm v}}}

\newcommand{\Ciso}{C_{\text{\rm iso}}}
\newcommand{\Cvol}{C_{\text{\rm vol}}}

\title[Random walk among random conductances]
{\large Functional CLT for random walk among\\bounded random conductances}
\author[M.~Biskup and T.~Prescott]{Marek Biskup and Timothy M. Prescott}

\begin{document}
\thanks{\hglue-4.5mm\fontsize{9.6}{9.6}\selectfont\copyright\,2007 by M.~Biskup and T.M.~Prescott. Reproduction, by any means, of the entire
article for non-commercial purposes is permitted without charge.\vspace{2mm}}
\maketitle

\vspace{-5mm}
\centerline{\textit{Department of Mathematics, University of California at Los Angeles}}

\vspace{-2mm}
\begin{abstract}
We consider the nearest-neighbor simple random walk on~$\Z^d$, $d\ge2$, driven by a field of i.i.d.\ random nearest-neighbor conductances $\omega_{xy}\in[0,1]$. Apart from the requirement that the bonds with positive conductances percolate, we pose no restriction on the law of the~$\omega$'s. We prove that, for a.e.\ realization of the environment, the path distribution of the walk converges weakly to that of non-degenerate, isotropic Brownian motion. The quenched functional CLT holds despite the fact that the local CLT may fail in~$d\ge5$ due to anomalously slow decay of the probability that the walk returns to the starting point at a given time.
\end{abstract}

\section{Introduction}
\label{sec1}\noindent
Let~$\B_d$ denote the set of unordered nearest-neighbor pairs (i.e., edges) of~$\Z^d$ and let $(\omega_b)_{b\in\B_d}$ be i.i.d.\ random variables with~$\omega_b\in[0,1]$. We will refer to $\omega_b$ as the \emph{conductance} of the edge~$b$. Let~$\BbbP$ denote the law of the~$\omega$'s and suppose that
\begin{equation}
\label{supercritical}
\BbbP(\omega_b>0)>p_\cc(d),
\end{equation}
where~$p_\cc(d)$ is the threshold for bond percolation on~$\Z^d$; in~$d=1$ we have~$p_\cc(d)=1$ so there we suppose~$\omega_b>0$ a.s. This condition ensures the existence of a unique infinite connected component~$\scrC_\infty$ of edges with strictly positive conductances; we will typically restrict attention to $\omega$'s for which~$\scrC_\infty$ contains a given site (e.g., the origin).

Each realization of $\scrC_\infty$ can be used to define a random walk~$X=(X_n)$ which moves about~$\scrC_\infty$ by picking, at each unit time, one of its $2d$ neighbors at random and passing to it with probability equal to the conductance of the corresponding edge. Technically, $X$ is a Markov chain with state space~$\scrC_\infty$ and the transition probabilities defined by
\begin{equation}
\label{p(x,y)}
P_{\omega,z}(X_{n+1}=y|X_n=x):=\frac{\omega_{xy}}{2d}
\end{equation}
if~$x,y\in\scrC_\infty$ and~$|x-y|=1$, and
\begin{equation}
\label{p(x,x)}
P_{\omega,z}(X_{n+1}=x|X_n=x):=1-\frac1{2d}\sum_{y\colon|y-x|=1}\omega_{xy}.
\end{equation}
The second index on~$P_{\omega,z}$ marks the initial position of the walk, i.e.,
\begin{equation}
%\label{}
P_{\omega,z}(X_0=z):=1.
\end{equation}
The counting measure on~$\scrC_\infty$ is invariant and reversible for this Markov chain.

\smallskip
The $d=1$ walk is a simple, but instructive, exercise in harmonic analysis of reversible random walks in random environments. Let us quickly sketch the proof of the fact that, for a.e.\ $\omega$ sampled from a translation-invariant, ergodic law on~$(0,1]^{\B_d}$ satisfying
\begin{equation}
\label{2nd}
\E\Bigl(\frac1{\omega_b}\Bigr)<\infty,
\end{equation}
the walk scales to Brownian motion under the usual diffusive scaling of space and time. (Here and henceforth~$\E$ denotes expectation with respect to the environment distribution.) The derivation works even for unbounded conductances provided \twoeqref{p(x,y)}{p(x,x)} are modified accordingly.

Abbreviate~$C:=\E(1/\omega_b)$. The key step of the proof is to realize that
\begin{equation}
\label{h-coordinate}
\varphi_\omega(x):=x+\frac1C\sum_{n=0}^{x-1}\,\Bigl(\frac1{\omega_{n,n+1}}-C\Bigr)
\end{equation}
is harmonic for the Markov chain. Hence $\varphi_\omega(X_n)$ is a martingale whose increments are, by \eqref{2nd} and a simple calculation, square integrable in the sense
\begin{equation}
\label{1d-diff}
\E E_{\omega,0}\bigl[\varphi_\omega(X_1)^2\bigr]<\infty.
\end{equation}
Invoking the stationarity and ergodicity of the Markov chain on the space of environments ``from the point of view of the particle'' --- we will discuss the specifics of this argument later --- the martingale $(\varphi_\omega(X_n))$ satisfies the conditions of the Lindeberg-Feller martingale functional CLT and so the law of $t\mapsto\varphi_\omega(X_{\lfloor nt\rfloor})/\sqrt n$ tends weakly to that of a Brownian motion with diffusion constant given by \eqref{1d-diff}. By the Pointwise Ergodic Theorem and \eqref{2nd} we have $\varphi_\omega(x)-x=o(x)$ as~$|x|\to\infty$. Thus the path~$t\mapsto X_{\lfloor nt\rfloor}/\sqrt n$ scales, in the limit~$n\to\infty$, to the same function as the deformed path $t\mapsto\varphi_\omega(X_{\lfloor nt\rfloor})/\sqrt n$. As this holds for a.e.~$\omega$, we have proved a \emph{quenched} functional CLT.

\smallskip
While the main ideas of the above $d=1$ solution work in all dimensions, the situation in~$d\ge2$ is, even for i.i.d.\ conductances, significantly more complicated. Progress has been made under additional conditions on the environment law. One such condition is \emph{strong ellipticity}, 
\begin{equation}
%\label{}
\exists\alpha>0\,\colon\quad\BbbP(\alpha\le\omega_b\le\ffrac1\alpha)=1.
\end{equation}
Here an annealed invariance principle was proved by Kipnis and Varadhan~\cite{Kipnis-Varadhan} and its queneched counterpart by Sidoravicius and Sznitman~\cite{Sidoravicius-Sznitman}. Another natural family of environments are those arising from supercritical \emph{bond percolation} on~$\Z^d$ for which~$(\omega_b)$ are i.i.d.\ zero-one valued with~$\BbbP(\omega_b=1)>p_\cc(d)$. For these cases an annealed invariance principle was proved by De Masi, Ferrari, Goldstein and Wick~\cite{demas1,demas2} and the quenched case was established in~$d\ge4$ by Sidoravicius and Sznitman~\cite{Sidoravicius-Sznitman}, and in all~$d\ge2$ by Berger and Biskup~\cite{BB} and Mathieu and Piatnitski~\cite{Mathieu-Piatnitski} . 

A common feature of the latter proofs is that, in~$d\ge3$, they require the use of heat-kernel upper bounds of the form
\begin{equation}
\label{heat-kernel}
P_{\omega,x}(X_n=y)\le \frac{c_1}{n^{d/2}}\exp\Bigl\{-c_2\frac{|x-y|^2}n\Bigr\},
\qquad x,y\in\scrC_\infty,
\end{equation}
where~$c_1,c_2$ are absolute constants and~$n$ is assumed to exceed a random quantity depending on the environment in the vicinity of~$x$ and~$y$. These were deduced by Barlow~\cite{Barlow} using sophisticated arguments that involve isoperimetry, regular volume growth and comparison of graph-theoretical and Euclidean distances for the percolation cluster. While the use of \eqref{heat-kernel} is conceptually rather unsatisfactory --- one seems to need a local-CLT level of control to establish a plain CLT --- no arguments (in~$d\ge3$) that avoid heat-kernel estimates are known at present.

The reliance on heat-kernel bounds also suffers from another problem: \eqref{heat-kernel} may actually \emph{fail} once the conductance law has sufficiently heavy tails at zero. This was first noted to happen by Fontes and Mathieu~\cite{Fontes-Mathieu} for the heat-kernel averaged over the environment; the quenched situation was analyzed recently by Berger, Biskup, Hoffman and Kozma~\cite{BBHK}. The main conclusion of~\cite{BBHK} is that the diagonal (i.e.,~$x=y$) bound in \eqref{heat-kernel} holds in~$d=2,3$ but the decay can be slower than any~$o(n^{-2})$ sequence in $d\ge5$. (The threshold sequence in~$d=4$ is presumably~$o(n^{-2}\log n)$.) This is caused by the existence of \emph{traps} that may capture the walk for a long time and thus, paradoxically, increase its chances to arrive back to the starting point.

The aformentioned facts lead to a natural question: In the absence of heat-kernel estimates, does the quenched CLT still hold? Our answer to this question is affirmative and constitutes the main result of this note. Another interesting question is what happens when the conductances are unbounded from above; this is currently being studied by Barlow and~Deuschel~\cite{Barlow-Deuschel}.

\smallskip\noindent
\emph{Note}: While this paper was in the process of writing, we received a preprint from Pierre Mathieu~\cite{Mathieu-CLT} in which he proves a result that is a continuous-time version of our main theorem. 
The strategy of~\cite{Mathieu-CLT} differs from ours by the consideration of a time-changed process (which we use only marginally) and proving that the ``new'' and ``old'' time scales are commensurate. Our approach is focused on proving the (pointwise) sublinearity of the corrector and it streamlines considerably the proof of~\cite{BB} in $d\ge3$ in that it limits the use of ``heat-kernel technology'' to a uniform bound on the heat-kernel decay (implied by isoperimetry) and a diffusive bound on the expected distance travelled by the walk (implied by regular volume growth).

\section{Main results and outline}
\noindent
Let~$\Omega:=[0,1]^{\B_d}$ be the set of all admissible random environments and let~$\BbbP$ be an i.i.d.\ law on~$\Omega$. Assuming \eqref{supercritical}, let $\scrC_\infty$ denote the a.s.\ unique infinite connected component of edges with positive conductances and introduce the conditional measure
\begin{equation}
%\label{}
\BbbP_0(-):=\BbbP(-|0\in\scrC_\infty).
\end{equation}
For~$T>0$, let $(C[0,T],\scrW_T)$ be the space of continuous functions~$f\colon[0,T]\to\R^d$ equipped with the Borel $\sigma$-algebra defined relative to the supremum topology. 

\smallskip
Here is our main result:

\begin{theorem}
\label{thm-main}
Suppose $d\ge2$ and $\BbbP(\omega_b>0)>p_\cc(d)$. For~$\omega\in\{0\in\scrC_\infty\}$, let~$(X_n)_{n\ge0}$ be the random walk with law~$P_{\omega,0}$ and let
\begin{equation}
\label{1.5a}
B_n(t):=\frac1{\sqrt n}\bigl(
X_{\lfloor tn\rfloor}+(tn-\lfloor tn\rfloor)(X_{\lfloor tn\rfloor+1}-X_{\lfloor tn\rfloor})\bigr),
\qquad t\ge0.
\end{equation}
Then for all~$T>0$ and for~$\BbbP_0$-almost every~$\omega$, the law of~$(B_n(t)\colon 0\le t\le T)$ on $(C[0,T],\scrW_T)$ converges, as~$n\to\infty$, weakly to the law of an isotropic Brownian motion $(B_t\colon 0\le t\le T)$ with a positive and finite diffusion constant (which is independent of~$\omega$).
\end{theorem}

Using a variant of~\cite[Lemma~6.4]{BB}, from here we can extract a corresponding conclusion for the ``agile'' version of our random walk (cf~\cite[Theorem~1.2]{BB}) by which we mean the walk that jumps from $x$ to its neighbor~$y$ with probability~$\omega_{xy}/\pi_\omega(x)$ where $\pi_\omega(x)$ is the sum of~$\omega_{xz}$ over all of the neighbors~$z$ of~$x$. Replacing discrete times by sums of i.i.d.\ exponential random variables, these invariance principles then extend also to the corresponding continuous-time processes. Theorem~\ref{thm-main} of course implies also an annealed invariance principle, which is the above convergence for the walk sampled from the path measure integrated over the environment.

\begin{remark}
As we were reminded by Y.~Peres, the above functional CLT automatically implies the ``usual'' lower bound on the heat-kernel. Indeed, the Markov property and reversibility of~$X$ yield
\begin{equation}
P_{\omega,0}(X_{2n}=0)\ge\sum_{\begin{subarray}{c}
x\in\scrC_\infty\\|x|\le\sqrt n
\end{subarray}}
P_{\omega,0}(X_n=x)^2.
\end{equation}
Cauchy-Schwarz then gives
\begin{equation}
%\label{}
P_{\omega,0}(X_{2n}=0)\ge P_{\omega,0}\bigl(|X_n|\le\sqrt n\bigr)^2\frac1{\bigl|\scrC_\infty\cap[-\sqrt n,\sqrt n]^d\bigr|}.
\end{equation}
Now Theorem~\ref{thm-main} implies that $P_{\omega,0}(|X_n|\le\sqrt n)$ is uniformly positive as~$n\to\infty$ and the Spatial Ergodic Theorem shows that $|\scrC_\infty\cap[-\sqrt n,\sqrt n]^d|$ grows proportionally to~$n^{d/2}$. Hence we get
\begin{equation}
%\label{}
P_{\omega,0}(X_{2n}=0)\ge\frac{C(\omega)}{n^{d/2}},\qquad n\ge1,
\end{equation}
with~$C(\omega)>0$ a.s. on the set~$\{0\in\scrC_\infty\}$. Note that, in $d=2,3$, this complements nicely the ``universal'' upper bounds derived in~\cite{BBHK}.
\end{remark}

\smallskip
The remainder of this paper is devoted to the proof of Theorem~\ref{thm-main}. The main line of attack is similar to the above 1D solution: We define a harmonic coordinate~$\varphi_\omega$ --- an analogue of \eqref{h-coordinate} --- and then prove an a.s.\ invariance principle for
\begin{equation}
%\label{}
t\mapsto\varphi_\omega(X_{\lfloor nt\rfloor})/\sqrt n
\end{equation}
along the martingale argument sketched before. The difficulty comes with showing the sublinearity of the corrector,
\begin{equation}
%\label{}
\varphi_\omega(x)-x=o(x),\qquad |x|\to\infty.
\end{equation}
As in Berger and Biskup~\cite{BB}, sublinearity can be proved directly along coordinate directions by soft ergodic-theory arguments. The crux is to extend this to a bound throughout $d$-dimensional boxes. 

Following the~$d\ge3$ proof of~\cite{BB}, the bound along coordinate axes readily implies \emph{sublinearity on average}, meaning that the set of~$x$ where $|\varphi_\omega(x)-x|$ exceeds $\epsilon|x|$ has zero density. The extension of sublinearity on average to pointwise sublinearity is the main novel part of the proof which, unfortunately, still makes non-trivial use of the ``heat-kernel technology.'' A heat-kernel upper bound of the form \eqref{heat-kernel} would do but, to minimize the extraneous input, we show that it suffices to have a diffusive bound for the expected displacement of the walk from its starting position. This step still requires detailed control of isoperimetry and volume growth as well as the comparison of the graph-theoretic and Euclidean distances, but it avoids many spurious calculations that are needed for the full-fledged heat-kernel estimates.

Of course, the required isoperimetric inequalities may not be true on~$\scrC_\infty$ because of the presence of ``weak'' bonds. As in~\cite{BBHK} we circumvent this by observing the random walk on the set of sites that have a connection to infinity by bonds with \emph{uniformly} positive conductances. Specifically we pick~$\alpha>0$ and let~$\scrC_{\infty,\alpha}$ denote the set of sites in~$\Z^d$ that are connected to infinity by a path whose edges obey $\omega_b\ge\alpha$. Here we note:

\begin{proposition}
\label{prop-percolace}
Let $d\ge2$ and~$p=\BbbP(\omega_b>0)>p_\cc(d)$. Then there exists $c(p,d)>0$ such that if~$\alpha$ satisfies
\begin{equation}
\label{infinite}
\BbbP(\omega_b\ge\alpha)>p_\cc(d)
\end{equation}
and
\begin{equation}
\label{finite}
\BbbP(0<\omega_b<\alpha)<c(p,d)
\end{equation}
then $\scrC_{\infty,\alpha}$ is nonempty and $\scrC_\infty\setminus\scrC_{\infty,\alpha}$ has only finite components a.s. In fact, if $\scrF(x)$ is the set of sites (possibly  empty) in the finite component of $\scrC_\infty\setminus\scrC_{\infty,\alpha}$ containing~$x$, then
\begin{equation}
\label{exponential}
\BbbP\bigl(x\in\scrC_\infty\,\,\&\,\,\diam \scrF(x)\ge n\bigr)\le C\texte^{-\eta n},\qquad n\ge1,
\end{equation}
for some~$C<\infty$ and $\eta>0$. Here ``$\diam$'' is the diameter in the~$\ell_\infty$ distance on~$\Z^d$.
\end{proposition}

The restriction of $\varphi_\omega$ to~$\scrC_{\infty,\alpha}$ is still harmonic, but with respect to a walk that can ``jump the holes'' of~$\scrC_{\infty,\alpha}$. A discrete-time version of this walk was utilized heavily in~\cite{BBHK}; for the purposes of this paper it will be more convenient to work with its continuous-time counterpart $Y=(Y_t)_{t\ge0}$. Explicitly, sample a path of the random walk~$X=(X_n)$ from~$P_{\omega,0}$ and denote by~$T_1,T_2,\dots$ the time intervals between successive visits of $X$ to~$\scrC_{\infty,\alpha}$. These are defined recursively by
\begin{equation}
\label{T-formula}
T_{j+1}:=\inf\,\bigl\{n\ge1\colon X_{T_1+\cdots+T_j+n}\in\scrC_{\infty,\alpha}\bigr\},
\end{equation}
with~$T_0=0$. For each~$x,y\in\scrC_{\infty,\alpha}$, let
\begin{equation}
\hat\omega_{xy}:=P_{\omega,x}(X_{T_1}=y)
\end{equation}
and define the operator
\begin{equation}
\label{generator}
(\LL_\omega^{(\alpha)} f)(x):=\sum_{y\in\scrC_{\infty,\alpha}}\hat\omega_{xy}\,\bigl[f(y)-f(x)\bigr].
\end{equation}
The continuous-time random walk~$Y$ is a Markov process with this generator; alternatively take the standard Poisson process $(N_t)_{t\ge0}$ with jump-rate one and set
\begin{equation}
Y_t:=X_{T_1+\cdots+T_{N_t}}.
\end{equation}
Note that, while~$Y$ may jump ``over the holes'' of~$\scrC_{\infty,\alpha}$, Proposition~\ref{prop-percolace} ensures that all of its jumps are finite. The counting measure on~$\scrC_{\infty,\alpha}$ is still invariant for this random walk, $\LL_\omega^{(\alpha)}$ is self-adjoint on the corresponding space of square integrable functions and~$\LL_\omega^{(\alpha)}\varphi_\omega=0$ on~$\scrC_{\infty,\alpha}$ (see Lemma~\ref{lemma-harmonic}).

\smallskip
The skeleton of the proof is condensed into the following statement whose proof, and adaptation to the present situation, is the main novel part of this note:

\begin{theorem}
\label{thm-sublinear}
Fix~$\alpha$ as in \twoeqref{infinite}{finite} and suppose~$\psi_\omega\colon\scrC_{\infty,\alpha}\to\R^d$ is a function and $\theta>0$ is a number such that the following holds for a.e.~$\omega$:
\settowidth{\leftmargini}{(111a)}
\begin{enumerate}
\item[(1)] (Harmonicity)
If~$\varphi_\omega(x):=x+\psi_\omega(x)$, then $\LL_\omega^{(\alpha)}\varphi_\omega=0$ on $\scrC_{\infty,\alpha}$.
\item[(2)] (Sublinearity on average)
For every~$\epsilon>0$,
\begin{equation}
\label{on-average}
\lim_{n\to\infty}\frac1{n^d}\sum_{\begin{subarray}{c}
x\in\scrC_{\infty,\alpha}\\|x|\le n
\end{subarray}}
\1_{\{|\psi_\omega(x)|\ge\epsilon n\}}=0.
\end{equation}
\item[(3)] (Polynomial growth)
\begin{equation}
\label{polygrowth}
\lim_{n\to\infty}\,\max_{\begin{subarray}{c}
x\in\scrC_{\infty,\alpha}\\|x|\le n
\end{subarray}}
\frac{|\psi_\omega(x)|}{n^\theta}=0.
\end{equation}
\end{enumerate}
Let~$Y=(Y_t)$ be the continuous-time random walk on~$\scrC_{\infty,\alpha}$ with generator~$\LL_\omega^{(\alpha)}$ and suppose~also:
\begin{enumerate}
\item[(4)] (Diffusive upper bounds) For a deterministic sequence~$b_n=o(n^2)$ and a.e.~$\omega$,
\begin{equation}
\label{diffusive}
\sup_{n\ge1}\,\,\max_{\begin{subarray}{c}
x\in\scrC_{\infty,\alpha}\\|x|\le n
\end{subarray}}\,\,
\sup_{t\ge b_n}\,\frac{E_{\omega,x}|Y_t-x|}{\sqrt t}<\infty
\end{equation}
and
\begin{equation}
\label{on-diag}
\sup_{n\ge1}\,\,\max_{\begin{subarray}{c}
x\in\scrC_{\infty,\alpha}\\|x|\le n
\end{subarray}}\,\,
\sup_{t\ge b_n}\, t^{d/2}P_{\omega,x}(Y_t=x)<\infty.
\end{equation}
\end{enumerate}
Then for almost every~$\omega$,
\begin{equation}
\label{sublinear}
\lim_{n\to\infty}\,\max_{\begin{subarray}{c}
x\in\scrC_{\infty,\alpha}\\|x|\le n
\end{subarray}}
\frac{|\psi_\omega(x)|}n=0.
\end{equation}
\end{theorem}

This result --- with $\psi_\omega$ playing the role of the corrector --- shows that $\varphi_\omega(x)-x=o(x)$ on~$\scrC_{\infty,\alpha}$. This readily extends to sublinearity on~$\scrC_\infty$ by the maximum principle applied to~$\varphi_\omega$ on the finite components of~$\scrC_\infty\setminus\scrC_{\infty,\alpha}$ and using that the component sizes obey a polylogarithmic upper bound. The assumptions (1-3) are known to hold for the corrector of the supercritical bond-percolation cluster and the proof applies, with minor modifications, to the present case as well. The crux is to prove \twoeqref{diffusive}{on-diag} which is where we need to borrow ideas from the ``heat-kernel technology.'' For our purposes it will suffice to take~$b_n=n$ in part~(4).

We remark that the outline strategy of proof extends rather seamlessly to other (tran\-sla\-tion-invariant, ergodic) conductance distributions with conductances bounded from above. Of course, one has to assume a number of specific properties for the ``strong'' component~$\scrC_{\infty,\alpha}$ that, in the i.i.d.\ case, we are able to check explicitly. 

\smallskip
The plan of the rest of this paper is a follows: Sect.~\ref{sec2b} is devoted to some basic percolation estimates needed in the rest of the paper. In Sect.~\ref{sec3} we define and prove some properties of the corrector $\chi$, which is a random function marking the difference between the harmonic coordinate $\varphi_\omega(x)$ and the geometric coordinate~$x$. In Sect.~\ref{sec4} we establish the a.s.\ sublinearity of the corrector as stated in Theorem~\ref{thm-sublinear} subject to the diffusive bounds \twoeqref{diffusive}{on-diag}. Then we assemble all facts into the proof of Theorem~\ref{thm-main}. Finally, in Sect.~\ref{sec5} we adapt some arguments from Barlow~\cite{Barlow} to prove \twoeqref{diffusive}{on-diag}; first in rather general Propositions~\ref{lemma-upper} and~\ref{thm-distance} and then for the case at hand.

\section{Percolation estimates}
\label{sec2b}\noindent
In this section we provide a proof of Proposition~\ref{prop-percolace} and also of a lemma dealing with the maximal distance the random walk~$Y$ can travel in a given number of steps. We will need to work with the ``static'' renormalization (cf~Grimmett \cite[Section~7.4]{Grimmett}) whose salient features we will now recall. The underlying ideas go back to the work of Kesten and Zhang~\cite{Kesten-Zhang}, Grimmett and Marstrand~\cite{Grimmett-Marstrand} and Antal and Pisztora~\cite{Antal-Pisztora}.

\smallskip
We say that an edge~$b$ is occupied if~$\omega_b>0$. Consider the lattice cubes
\begin{equation}
%\label{}
B_L(x):=x+[0,L]^d\cap\Z^d
\end{equation}
and
\begin{equation}
\tilde B_{3L}(x):=x+[-L,2L]^d\cap\Z^d
\end{equation}
and note that~$\tilde B_{3L}(x)$ consists of $3^d$ copies of~$B_L(x)$ that share only sites on their adjacent boundaries. Let~$G_L(x)$ be the ``good event'' --- whose occurrence designates~$B_L(Lx)$ to be a ``good block'' --- which is the set of configurations such that:
%\settowidth{\leftmargini}{(11)}
\begin{enumerate}
\item[(1)]
For each neighbor~$y$ of~$x$, the side of the block~$B_L(Ly)$ adjacent to~$B_L(Lx)$ is connected to the opposite side of~$B_L(Ly)$ by an occupied path.
\item[(2)]
Any two occupied paths connecting $B_L(Lx)$ to the boundary of $\tilde B_{3L}(Lx)$ are connected by an occupied path using only edges with both endpoints in~$\tilde B_{3L}(Lx)$.
\end{enumerate}
The sheer existence of infinite cluster implies that (1) occurs with high probability once~$L$ is large (see Grimmett~\cite[Theorem~8.97]{Grimmett}) while the situation in~(2) occurs with large probability once there is percolation in half space (see Grimmett~\cite[Lemma~7.89]{Grimmett}). It follows that
\begin{equation}
%\label{}
\BbbP\bigl(G_L(x)\bigr)\,\underset{L\to\infty}\longrightarrow\,1
\end{equation}
whenever~$\BbbP(\omega_b>0)>p_\cc(d)$. A crucial consequence of the above conditions is that, if~$G_L(x)$ and~$G_L(y)$ occur for neighboring sites~$x,y\in\Z^d$, then the largest connected components in~$\tilde B_{3L}(Lx)$ and~$\tilde B_{3L}(Ly)$ --- sometimes referred to as \emph{spanning clusters} --- are connected. Thus, if~$G_L(x)$ occurs for all~$x$ along an infinite path on~$\Z^d$, the corresponding spanning clusters are subgraphs of~$\scrC_\infty$.

A minor complication is that the events~$\{G_L(x)\colon x\in\Z^d\}$ are not independent. However, they are $4$-dependent in the sense that if $(x_i)$ and~$(y_j)$ are such that $|x_i-y_j|>4$ for each~$i$ and~$j$, then the families~$\{G_L(x_i)\}$ and~$\{G_L(y_j)\}$ are independent. By the main result of Liggett, Schonmann and Stacey~\cite[Theorem~0.0]{LSS} (cf~\cite[Theorem~7.65]{Grimmett}) the indicators~$\{\1_{G_L(x)}\colon x\in\Z^d\}$, regarded as a random process on~$\Z^d$, stochastically dominate i.i.d.\ Bernoulli random variables whose density (of ones) tends to one as~$L\to\infty$.

\begin{proofsect}{Proof of Proposition~\ref{prop-percolace}}
In~$d=2$ the proof is actually very simple because it suffices to choose $\alpha$ such that \eqref{infinite} holds. Then~$\scrC_\infty\setminus\scrC_{\infty,\alpha}\subset\Z^2\setminus\scrC_{\infty,\alpha}$ has only finite (subcritical) components whose diameter has exponential tails \eqref{exponential} by, e.g.,~\cite[Theorem~6.10]{Grimmett}. 

\newcommand{\pomega}{\partial^{\,\omega}\!}

To handle general dimensions we will have to invoke the above ``static'' renormalization. Let $G_L(x)$ be as above and consider the event~$G_{L,\alpha}(x)$ where we in addition require that~$\omega_b\not\in(0,\alpha)$ for every edge with both endpoints in~$\tilde B_{3L}(Lx)$. Clearly, 
\begin{equation}
\label{}
\lim_{L\to\infty}\lim_{\alpha\downarrow0}\BbbP\bigl(G_{L,\alpha}(x)\bigr)=1.
\end{equation}
Using the aforementioned domination by site percolation, and adjusting~$L$ and~$\alpha$ we can thus ensure that, with probability one, the set
\begin{equation}
\label{GL-set}
\bigl\{x\in\Z^d\colon G_{L,\alpha}(x)\text{ occurs}\bigr\}
\end{equation}
has a unique infinite component $\eusm C_\infty$, whose complement has only finite components. Moreover, if~$\eusm G(0)$ is the finite connected component of~$\Z^d\setminus\eusm C_\infty$ containing the origin, then a standard Peierls argument yields
\begin{equation}
\label{G-tail}
\BbbP\bigl(\diam\eusm G(0)\ge n\bigr)\le\texte^{-\zeta n}
\end{equation}
for some~$\zeta>0$. To prove \eqref{exponential}, it suffices to show that 
\begin{equation}
\label{finite-contained}
\scrF(0)\subset\bigcup_{x\in\eusm G(0)} B_L(Lx)
\end{equation}
once~$\diam\scrF(0)>3L$. Indeed, then $\diam\scrF(0)\le L\diam\eusm G(0)$ and so \eqref{G-tail} implies \eqref{exponential} with~$\eta:=\ffrac\zeta L$ and~$C:=\texte^{3L\eta}$.

To prove \eqref{finite-contained}, pick $z\in\scrF(0)$ and let $x$ be such that $z\in B_L(Lx)$. It suffices to show that if~$G_{L,\alpha}(x)$ occurs, then $x$ is not adjacent to an infinite component in \eqref{GL-set}. Assuming that~$x$ \emph{is} adjacent to such a component, the fact that the spanning clusters in adjecent ``good blocks'' are connected and thus contained in~$\scrC_{\infty,\alpha}$ implies
\begin{equation}
%\label{}
\scrC_{\infty,\alpha}\cap B_L(Lx)\ne\emptyset.
\end{equation}
But then~$B_L(Lx)$ is intersected by two ``large'' components, $\scrC_{\infty,\alpha}$ and~$\scrF(0)$, of edges with~$\omega_b\ge\alpha$. (This is where we need $\diam\scrF(0)>3L$.) If these components are joined by an occupied path --- i.e., a path of edges with~$\omega_b>0$ --- within~$\tilde B_{3L}(Lx)$, then~$\tilde B_{3L}(Lx)$ contains a ``weak'' bond and so~$G_{L,\alpha}$ fails. In the absence of such a path the requirement~(2) in the definition of~$G_L(x)$ is not satisfied and so~$G_{L,\alpha}(x)$ fails too.
\end{proofsect}

Let~$\textd(x,y)$ be the ``Markov distance'' on~$V=\scrC_{\infty,\alpha}$, i.e., the minimal number of jumps the random walk~$Y=(Y_t)$ needs to make to get from~$x$ to~$y$. Note that~$\textd(x,y)$ could be quite smaller than the graph-theoretic distance on~$\scrC_{\infty,\alpha}$ and/or the Euclidean distance. (The latter distances are known to be comparable, see Antal and Pisztora~\cite{Antal-Pisztora}.) To control the volume growth for the Markov graph of the random walk~$Y$ we will need to know that~$\textd(x,y)$ is nevertheless comparable with the Euclidean distance~$|x-y|$:

\begin{lemma}
\label{lemma-compare}
There exists~$\varrho>0$ and for each~$\gamma>0$ there is $\alpha>0$ obeying \twoeqref{infinite}{finite} and~$C<\infty$ such that
\begin{equation}
%\label{}
\BbbP\bigl(\,0,x\in\scrC_{\infty,\alpha}\,\,\&\,\,\textd(0,x)\le \varrho|x|\,\bigr)\le C\texte^{-\gamma|x|},
\qquad x\in\Z^d.
\end{equation}
\end{lemma}

\begin{proofsect}{Proof}
Suppose~$\alpha$ is as in the proof of Proposition~\ref{prop-percolace}. Let~$(\eta_x)$ be independent Bernoulli that dominate the indicators $\1_{G_{L,\alpha}}$ from below and let~$\eusm C_\infty$ be the unique infinite component of the set $\{x\in\Z^d\colon\eta_x=1\}$. We may ``wire'' the ``holes'' of~$\eusm C_\infty$ by putting an edge between every pair of sites on the external boundary of each finite component of~$\Z^d\setminus\eusm C_\infty$; we use~$\textd'(0,x)$ to denote the distance between~$0$ and~$x$ on the induced graph. The processes~$\eta$ and $(\1_{G_{L,\alpha}(x)})$ can be coupled so that each connected component of~$\scrC_\infty\setminus\scrC_{\infty,\alpha}$ with diameter exceeding~$3L$ is ``covered'' by a finite component of~$\Z^d\setminus\eusm C_\infty$, cf~\eqref{finite-contained}. As is easy to check, this implies
\begin{equation}
\label{dd-prime}
\textd(0,x)\ge\textd'(0,x')
\quad\text{and}\quad
|x'|\ge\frac1L|x|-1
\end{equation}
whenever~$x\in B_L(Lx')$. It thus suffices to show the above bound for distance~$\textd'(0,x')$.

Let~$p=p_{L,\alpha}$ be the parameter of the Bernoulli distribution and recall that~$p$ can be made as close to one as desired by adjusting~$L$ and~$\alpha$. Let~$z_0=0,z_1,\dots,z_n=x$ be a nearest-neighbor path on~$\Z^d$. Let~$\eusm G(z_i)$ be the unique finite component of~$\Z^d\setminus\eusm C_\infty$ that contains~$z_i$ --- if~$z_i\in\eusm C_\infty$, we have~$\eusm G(z_i)=\emptyset$. Define
\begin{equation}
%\label{}
\ell(z_0,\dots,z_n)\,:= \sum_{i=0}^n\diam \eusm G(z_i)\Bigl(\,\prod_{j<i}\1_{\{z_j\not\in \eusm G(z_i)\}}\Bigr).
\end{equation}
We claim that for each $\lambda>0$ we can adjust~$p$ so that
\begin{equation}
\label{exp-moment}
\E\,\texte^{\lambda\ell(z_0,\dots,z_n)}\le\texte^n
\end{equation}
for all~$n\ge1$ and all paths as above. To verify this we note that the components contributing to $\ell(z_0,\dots,z_n)$ are distance at least one from one another. So conditioning on all but the last component, and the sites in the ultimate vicinity, we may use the Peierls argument to estimate the conditional expectation of $\texte^{\lambda\diam \eusm G(z_n)}$ by, say,~$\texte^1$. (We are using also that~$\diam \eusm G(z_n)$ is smaller than the boundary of~$\eusm G(z_n)$.) Proceeding by induction~$n$ times, \eqref{exp-moment} follows.

For any given~$\gamma>0$ we can adjust~$p$ so that \eqref{exp-moment} holds for
\begin{equation}
%\label{}
\lambda:=2(1+\log(2d)+\gamma)
\end{equation}
 As the number of nearest-neighbor paths~$(z_0=0,\dots,z_n=x)$ is bounded by $(2d)^n$, an exponential Chebyshev estimate then shows
\begin{equation}
%\label{}
\BbbP\Bigl(\exists(z_0=0,\dots,z_n=x)\colon\ell(z_0,\dots,z_n)>\frac n2\Bigr)\le \texte^{-\gamma n}.
\end{equation}
But if~$(z_0=0,\dots,z_n=x)$ is the shortest nearest-neighbor interpolation of a path that achieves $\textd'(0,x)$, then
\begin{equation}
%\label{}
\textd'(0,x)\ge n-\ell(z_0,\dots,z_n).
\end{equation}
Since, trivially, $|x|\le n$, we deduce
\begin{equation}
%\label{}
\BbbP\bigl(\textd'(0,x)\le\tfrac12|x|\bigr)\le \texte^{-\gamma |x|}
\end{equation}
as desired.
\end{proofsect}

\section{Corrector}
\label{sec3}\noindent
The purpose of this section is to define, and prove some properties of, the \emph{corrector} $\chi(\omega,x):=\varphi_\omega(x)-x$. This object could be defined probabilistically by the limit
\begin{equation}
\chi(\omega,x)=\lim_{n\to\infty}\bigl(E_{\omega,x}(X_n)-E_{\omega,0}(X_n)\bigr)-x,
\end{equation}
unfortunately, at this moment we seem to have no direct (probabilistic) argument showing that the limit exists. The traditional definition of the corrector involves spectral calculus (Kipnis and Varadhan~\cite{Kipnis-Varadhan}); we will invoke a projection construction from Mathieu and Piatnitski~\cite{Mathieu-Piatnitski} (see also Giacomin, Olla and Spohn~\cite{Giacomin-Olla-Spohn}). 

Let~$\BbbP$ be an i.i.d.\ law on $(\Omega,\scrF)$ where~$\Omega:=[0,1]^{\B_d}$ and $\scrF$ is the natural product $\sigma$-algebra. Let~$\tau_x\colon\Omega\to\Omega$ denote the shift by~$x$, i.e.,
\begin{equation}
%\label{}
(\tau_z\omega)_{xy}:=\omega_{x+z,y+z}
\end{equation}
and note that~$\BbbP\circ\tau_x^{-1}=\BbbP$ for all~$x\in\Z^d$. Recall that~$\scrC_\infty$ is the infinite connected component of edges with~$\omega_b>0$ and, for $\alpha>0$, let~$\scrC_{\infty,\alpha}$ denote the set of sites connected to infinity by edges with $\omega_b\ge\alpha$. If $\BbbP(0\in\scrC_{\infty,\alpha})>0$, let
\begin{equation}
\BbbP_\alpha(-):=\BbbP(-|0\in\scrC_{\infty,\alpha})
\end{equation}
and let~$\E_\alpha$ be the corresponding expectation.  Given~$\omega\in\Omega$ and sites $x,y\in\scrC_{\infty,\alpha}(\omega)$, let~$\textd_\omega^{(\alpha)}(x,y)$ denote the graph distance between~$x$ and~$y$ as measured on~$\scrC_{\infty,\alpha}$. (Note this is distinct from the Markov distance $\textd(x,y)$ discussed, e.g., in Lemma~\ref{lemma-compare}.) We will also use $\LL_\omega$ to denote the generator of the continuous-time version of the walk~$X$, i.e.,
\begin{equation}
%\label{}
(\LL_\omega f)(x):=\frac1{2d}\sum_{y\colon|y-x|=1}\omega_{xy}\,\bigl[f(y)-f(x)\bigr].
\end{equation}
The following theorem summarizes all relevant properties of the corrector:

\begin{theorem}
\label{thm3.1}
Suppose~$\BbbP(0\in\scrC_\infty)>0$. There exists a function~$\chi\colon\Omega\times\Z^d\to\R^d$ such that the following holds $\BbbP_0$-a.s.:
\settowidth{\leftmargini}{(111a)}
\begin{enumerate}
\item[(1)] (Gradient field)
$\chi(0,\omega)=0$ and, for all $x,y\in\scrC_{\infty}(\omega)$,
\begin{equation}
\label{2.12d}
\chi(\omega,x)-\chi(\omega,y)=\chi(\tau_y\omega,x-y).
\end{equation}
\item[(2)] (Harmonicity)
$\varphi_\omega(x)\,:=x+\chi(\omega,x)$ obeys~$\LL_\omega\varphi_\omega=0$.
\item[(3)] (Square integrability)
There is a constant~$C=C(\alpha)<\infty$ such that for all $x,y\in\Z^d$ satisfying~$|x-y|=1$,
\begin{equation}
\label{eq:l2bnd}
\E_\alpha\bigl(\vert\chi(\cdot,y)-\chi(\cdot,x)\vert^2\,\omega_{xy}\1_{\{x\in\scrC_\infty\}}\bigr)<C
\end{equation}
\end{enumerate}
Let~$\alpha>0$ be such that $\BbbP(0\in\scrC_{\infty,\alpha})>0$. Then we also have:
\begin{enumerate}
\item[(4)] (Polynomial growth)
For every $\theta>d$, a.s.,
\begin{equation}
%\label{}
\lim_{n\to\infty}\,\max_{\begin{subarray}{c}
x\in\scrC_{\infty,\alpha}\\|x|\le n
\end{subarray}}
\frac{|\chi(\omega,x)|}{n^\theta}=0.
\end{equation}
\item[(5)] (Zero mean under random shifts)
Let~$Z\colon\Omega\to\Z^d$ be a random variable such~that
\settowidth{\leftmarginii}{(1111)}
\begin{enumerate}
\item[(a)]
$Z(\omega)\in\scrC_{\infty,\alpha}(\omega)$,
\item[(b)]
$\BbbP_\alpha$ is preserved by~$\omega\mapsto\tau_{Z(\omega)}(\omega)$,
\item[(c)]
$\E_\alpha(\textd_\omega^{(\alpha)}(0,Z(\omega))^q)<\infty$ for some~$q>3d$.
\end{enumerate}
Then $\chi(\cdot,Z(\cdot))\in L^1(\Omega,\scrF,\BbbP_\alpha)$ and
\begin{equation}
%\label{}
\E_\alpha\bigl[\chi(\cdot,Z(\cdot))\bigr]=0.
\end{equation}
\end{enumerate}
\end{theorem}

As noted before, to construct the corrector we will invoke a projection argument.   Abbreviate $L^2(\Omega)=L^2(\Omega,\scrF,\BbbP_0)$ and let~$B:=\{\hate_1,\dots,\hate_d\}$ be the set of coordinate vectors.  Consider the space $L^2(\Omega\times B)$ of square integrable functions $u\colon\Omega\times B\to\R^d$ equipped with the inner product
\begin{equation}
%\label{}
(u,v):=\E_0\Bigl(\,\sum_{b\in B} u(\omega,b)\cdot v(\omega,b)\,\omega_b\Bigr).
\end{equation}
We may interpret $u\in L^2(\Omega\times B)$ as a flow by putting
$u(\omega,-b)=-u(\tau_{-b}\omega,b)$. 
Some, but not all, elements of~$L^2(\Omega\times B)$ can be obtained as gradients of local functions, where the \emph{gradient}~$\nabla$ is the map~$L^2(\Omega)\to L^2(\Omega\times B)$ defined by
\begin{equation}
%\label{}
(\nabla\phi)(\omega,b):=\phi\circ\tau_b(\omega)-\phi(\omega).
\end{equation}
Let~$L^2_\nabla$ denote the closure of the set of gradients of all \emph{local} functions --- i.e., those depending only on the portion of~$\omega$ in a finite subset of~$\Z^d$ --- and
note the following orthogonal decomposition~$L^2(\Omega\times B)=L^2_\nabla\oplus (L^2_\nabla)^\perp$.

The elements of~$(L^2_\nabla)^\perp$ can be characterized using the concept of \emph{divergence}, which for $u\colon\Omega\times B\to\R^d$ is the function~$\DIV u\colon\Omega\to\R^d$ defined by
\begin{equation}
%\label{}
\DIV u(\omega):=\sum_{b\in B}\bigl[\omega_bu(\omega,b)-\omega_{-b}u(\tau_{-b}\omega,b)\bigr].
\end{equation}
Using the interpretation of~$u$ as a flow, $\DIV u$ is simply the net flow out of the origin.
The characterization of $(L^2_\nabla)^\perp$ is now as follows:

\begin{lemma}
\label{lemma3.2}
$u\in (L^2_\nabla)^\perp$ if and only if $\DIV u(\omega)=0$ for $\BbbP_0$-a.e.~$\omega$.
\end{lemma}

\begin{proofsect}{Proof}
Let~$u\in L^2(\Omega\times B)$ and let $\phi\in L^2(\Omega)$ be a local function. A direct calculation and the fact that~$\omega_{-b}=(\tau_{-b}\omega)_{b}$ yield
\begin{equation}
%\label{}
(u,\nabla\phi)=-\E_0\Bigl(\,\phi(\omega)\,\DIV u(\omega)\Bigr).
\end{equation}
If~$u\in(L^2_\nabla)^\perp$, then $\DIV u$ integrates to zero against all local functions. Since these are dense in~$L^2(\Omega)$, we have~$\DIV u=0$ a.s.
\end{proofsect}

It is easy to check that every~$u\in L^2_\nabla$ is \emph{curl-free} in the sense that for any oriented loop $(x_0,x_1,\dots,x_n)$ on~$\scrC_\infty(\omega)$ with~$x_n=x_0$ we have
\begin{equation}
\label{cycle}
\sum_{j=0}^{n-1} u(\tau_{x_j}\omega,x_{j+1}-x_j)=0.
\end{equation}
On the other hand, every $u\colon\Omega\times B\to\R^d$ which is curl-free can be integrated into a unique function~$\phi\colon\Omega\times\scrC_\infty(\cdot)\to\R^d$ such that
\begin{equation}
%\label{}
\phi(\omega,x)=\sum_{j=0}^{n-1} u(\tau_{x_j}\omega,x_{j+1}-x_j)
\end{equation}
holds for any path~$(x_0,\dots,x_n)$ on~$\scrC_\infty(\omega)$ with $x_0=0$ and $x_n=x$. This function will automatically satisfy the \emph{shift-covariance} property
\begin{equation}
\label{shift-cov}
\phi(\omega,x)-\phi(\omega,y)=\phi(\tau_y\omega,x-y),\qquad x,y\in\scrC_\infty(\omega).
\end{equation}
We will denote the space of such functions~$\HH(\Omega\times\Z^d)$.
To denote the fact that~$\phi$ is assembled from the shifts of~$u$, we will write
\begin{equation}
%\label{}
u=\GRAD\phi,
\end{equation}
i.e.,``$\GRAD$'' is a map from $\HH(\Omega\times\Z^d)$ to functions $\Omega\times B\to\R^d$ that takes a function $\phi\in\HH(\Omega\times\Z^d)$ and assigns to it the collection of values $\{\phi(\cdot,b)-\phi(\cdot,0)\colon b\in B\}$.

\begin{lemma}
\label{lemma3.3}
Let~$\phi\in\HH(\Omega\times\Z^d)$ be such that~$\GRAD\phi\in (L^2_\nabla)^\perp$. Then~$\phi$ is (discrete) harmonic for the random walk on~$\scrC_\infty$, i.e., for~$\BbbP_0$-a.e.~$\omega$ and all~$x\in\scrC_\infty(\omega)$,
\begin{equation}
%\label{}
(\LL_\omega\phi)(\omega,x)=0.
\end{equation}
\end{lemma}

\begin{proofsect}{Proof}
Our definition of divergence is such that ``div grad = $2d\,\LL_\omega$'' holds.
Lemma~\ref{lemma3.2} implies that~$u\in(L^2_\nabla)^\perp$ if and only if~$\DIV u=0$, which is equivalent to $(\LL_\omega\phi)(\omega,0)=0$. By translation covariance, this extends to all sites in~$\scrC_\infty$.
\end{proofsect}

\begin{proofsect}{Proof of Theorem~\ref{thm3.1}(1-3)}
Consider the function~$\phi(\omega,x):=x$ and let~$u:=\GRAD\phi$. Clearly, $u\in L^2(\Omega\times B)$. Let $G\in L^2_\nabla$ be the orthogonal projection of~$-u$ onto~$L^2_\nabla$ and define $\chi\in\HH(\Omega\times\Z^d)$ to be the unique function such that
\begin{equation}
%\label{}
G=\GRAD\chi\quad\text{and}\quad \chi(\cdot,0)=0.
\end{equation}
This definition immediately implies \eqref{2.12d}, while the definition of the inner product on $L^2(\Omega\times B)$ directly yields~\eqref{eq:l2bnd}. Since~$u$ projects to~$-G$ on~$L^2_\nabla$, we have~$u+G\in(L^2_\nabla)^\perp$. But~$u+G=\GRAD[x+\chi(\omega,x)]$ and so, by Lemma~\ref{lemma3.3}, $x\mapsto x+\chi(\omega,x)$ is harmonic with respect to $\LL_\omega$.
\end{proofsect}

\begin{remark}
We note that the corrector is actually uniquely determined by properties (1-3) of Theorem~\ref{thm3.1}. In fact, $x+\chi$ spans the orthogonal complement of~$L^2_\nabla$ in the space of shift-covariant functions.
See Biskup and Spohn~\cite{Biskup-Spohn}.
\end{remark}

For the remaining parts of Theorem~\ref{thm3.1} we will need to work on~$\scrC_{\infty,\alpha}$. However, we do not yet need the full power of Proposition~\ref{prop-percolace}; it suffices to note that~$\scrC_{\infty,\alpha}$ has the law of a supercritical percolation cluster.

\begin{proofsect}{Proof of Theorem~\ref{thm3.1}(4)}
Let $\theta>d$ and abbreviate
\begin{equation}
\label{mn}
R_n:=\,\max_{\begin{subarray}{c}
x\in\scrC_{\infty,\alpha}\\|x|\le n
\end{subarray}}
\bigl|\chi(\omega,x)\bigr|.
\end{equation}
By Theorem~1.1 of Antal and Pisztora~\cite{Antal-Pisztora},
\begin{equation}
%\label{}
\lambda(\omega)\,:=\sup_{x\in\scrC_{\infty,\alpha}}\,\frac{\textd_\omega^{(\alpha)}(0,x)}{|x|}<\infty,
\qquad\BbbP_\alpha\text{\rm-a.s.},
\end{equation}
and so it suffices to show that $R_n/n^\theta\to0$ on $\{\lambda(\omega)\le\lambda\}$ for every~$\lambda<\infty$.
But on $\{\lambda(\omega)\le\lambda\}$ every~$x\in\scrC_{\infty,\alpha}$ with~$|x|\le n$ can be reached by a path on~$\scrC_{\infty,\alpha}$ that does not leave~$[-\lambda n,\lambda n]^d$ and so, on $\{\lambda(\omega)\le\lambda\}$,
\begin{equation}
%\label{}
R_n\le\sum_{\begin{subarray}{c}
x\in\scrC_{\infty,\alpha}\\|x|\le \lambda n
\end{subarray}}\sum_{b\in B}
\sqrt{\frac{\omega_{x,x+b}}\alpha}\bigl|\chi(\omega,x+b)-\chi(\omega,x)\bigr|.
\end{equation}
Invoking the bound \eqref{eq:l2bnd} we then get
\begin{equation}
%\label{}
\Vert R_n\1_{\{\lambda(\omega)\le\lambda\}}\Vert_2\le C n^d
\end{equation}
for some constant~$C=C(\alpha,\lambda,d)<\infty$. Applying Chebyshev's inequality and summing~$n$ over powers of~$2$ then yields $R_n/n^\theta\to0$ a.s.\ on $\{\lambda(\omega)\le\lambda\}$.
\end{proofsect}

\begin{proofsect}{Proof of Theorem~\ref{thm3.1}(5)}
Let~$Z$ be a random variable satisfying the properties~(a-c). By the fact that~$G\in L^2_\nabla$, there exists a sequence~$\psi_n\in L^2(\Omega)$ such that
\begin{equation}
%\label{}
\psi_n\circ\tau_x-\psi_n\,\underset{n\to\infty}\longrightarrow\,\chi(\cdot,x)\qquad\text{in }L^2(\Omega\times B).
\end{equation}
Abbreviate~$\chi_n(\omega,x)=\psi_n\circ\tau_x(\omega)-\psi_n(\omega)$ and without loss of generality assume that $\chi_n(\cdot,x)\to\chi(\cdot,x)$ almost surely.

By the fact that~$Z$ is $\BbbP_\alpha$-preserving we have~$\E_\alpha(\chi_n(\cdot,Z(\cdot)))=0$ as soon as we can show that~$\chi_n(\cdot,Z(\cdot))\in L^1(\Omega)$. It thus suffices to prove that
\begin{equation}
\label{L1conv}
\chi_n\bigl(\cdot,Z(\cdot)\bigr)\,\underset{n\to\infty}\longrightarrow\,\chi\bigl(\cdot,Z(\cdot)\bigr)
\qquad\text{in }L^1(\Omega).
\end{equation}
Abbreviate~$K(\omega):=\textd_\omega^{(\alpha)}(0,Z(\omega))$ and note that, as in part~(4),
\begin{equation}
%\label{}
\bigl|\chi_n(\omega,Z(\omega))\bigr|\le
\sum_{\begin{subarray}{c}
x\in\scrC_{\infty,\alpha}\\|x|\le K(\omega)
\end{subarray}}\sum_{b\in B}
\sqrt{\frac{\omega_{x,x+b}}\alpha}\bigl|\chi_n(\omega,x+b)-\chi_n(\omega,x)\bigr|.
\end{equation}
The quantities
\begin{equation}
%\label{}
\sqrt{\omega_{x,x+b}}\,\bigl|\,\chi_n(\omega,x+b)-\chi_n(\omega,x)\bigr|\1_{\{x\in\scrC_{\infty,\alpha}\}}
\end{equation}
are bounded in~$L^2$, uniformly in~$x$, $b$ and $n$, and assumption~(c) tells us that $K\in L^q$ for some~$q>3d$. Ordering the edges in~$\B_d$ according to their distance from the origin, Lemma~4.5 of Berger and Biskup~\cite{BB} with the specific choices
\begin{equation}
%\label{}
p:=2,\quad  s:=\ffrac qd\quad\text{and}\quad N:=d(2K+1)^d
\end{equation}
(note that~$N\in L^s(\Omega)$) implies that for some~$r>1$,
\begin{equation}
%\label{}
\sup_{n\ge1}\Vert\chi_n(\cdot,Z(\cdot))\Vert_r<\infty.
\end{equation}
Hence, the family $\{\chi_n(\cdot,Z(\cdot))\}$ is uniformly integrable and \eqref{L1conv} thus follows by the fact that~$\chi_n(\cdot,Z(\cdot))$ converge almost surely.
\end{proofsect}

\begin{remark}
It it worth pointing out that the proof of properties~(1-3) extends nearly verbatim to the setting with arbitrary conductances and arbitrary long jumps (i.e., the case when $B$ is simply all of~$\Z^d$). One only needs that~$x$ is in $L^2(\Omega\times B)$, i.e.,
\begin{equation}
%\label{}
\E\biggl(\,\sum_{x\in\Z^d}\omega_{0,x}|x|^2\biggr)<\infty.
\end{equation}
The proof of (4-5) seems to require additional (and somewhat unwieldy) conditions.
\end{remark}

\section{Convergence to Brownian motion}
\label{sec4}\noindent
Here we will prove Theorem~\ref{thm-main}. We commence by establishing the conclusion of Theorem~\ref{thm-sublinear} whose proof draws on an idea, suggested to us by Yuval Peres, that sublinearity on average plus heat kernel upper bounds imply pointwise sublinearity. We have reduced the extraneous input from heat-kernel technology to the assumptions \twoeqref{diffusive}{on-diag}. These imply heat-kernel upper bounds but generally require less work to prove.

\smallskip
The main technical part of Theorem~\ref{thm-main} is encapsulated into the following lemma:

\begin{lemma}
\label{lemma-encapsulate}
Abusing the notation from \eqref{mn} slightly, let
\begin{equation}
%\label{}
R_n:=\,\max_{\begin{subarray}{c}
x\in\scrC_{\infty,\alpha}\\|x|\le n
\end{subarray}}
\bigl|\psi_\omega(x)\bigr|.
\end{equation}
Under the conditions (1,2,4) of Theorem~\ref{thm-sublinear}, for each~$\epsilon>0$ and $\delta>0$, there exists an a.s.\ finite random variable~$n_0=n_0(\omega,\epsilon,\delta)$ such that
\begin{equation}
\label{mn-bd}
R_n\le\epsilon n+\delta R_{3n}.
\qquad n\ge n_0.
\end{equation}
\end{lemma}

Before we prove this, let us see how this and \eqref{polygrowth} imply \eqref{sublinear}.

\begin{proofsect}{Proof of Theorem~\ref{thm-sublinear}}
Suppose that $R_n/n\not\to0$ and pick~$c$ with $0<c<\limsup_{n\to\infty}R_n/n$. Let $\theta$ be is as in \eqref{polygrowth} and choose
\begin{equation}
%\label{}
\epsilon:=\frac c2\quad\text{and}\quad\delta:=\frac1{3^{\theta+1}}.
\end{equation}
Note that then $c'-\epsilon\ge3^\theta\delta c'$ for all~$c'\ge c$. If~$R_n\ge c n$ --- which happens for infinitely many~$n$'s --- and~$n\ge n_0$, then \eqref{mn-bd} implies
\begin{equation}
%\label{}
R_{3n}\ge\frac{c-\epsilon}{\delta}n\ge3^\theta c n
\end{equation}
and, inductively, $R_{3^kn}\ge 3^{k\theta}c n$. However, that contradicts \eqref{polygrowth} by which $R_{3^kn}/3^{k\theta}\to0$ as $k\to\infty$ (with~$n$ fixed).
\end{proofsect}

The idea underlying Lemma~\ref{lemma-encapsulate} is simple: We run a continuous-time random walk $(Y_t)$ for time $t=o(n^2)$ starting from the maximizer of~$R_n$ and apply the harmonicity of $x\mapsto x+\psi_\omega(x)$ to derive an estimate on the expectation of~$\psi(Y_t)$. The right-hand side of \eqref{mn-bd} expresses two characteristic situations that may occur at time~$t$: Either we have $|\psi_\omega(Y_t)|\le\epsilon n$ --- which, by ``sublinearity  on average,'' happens with overwhelming probability --- or $Y$ will not yet have left the box $[-3n,3n]^d$ and so $\psi_\omega(Y_t)\le R_{3n}$. The point is to show that these are the dominating strategies.

\begin{proofsect}{Proof of Lemma~\ref{lemma-encapsulate}}
Fix~$\epsilon,\delta>0$ and let~$C_1=C_1(\omega)$ and~$C_2=C_2(\omega)$ denote the suprema in \eqref{diffusive} and \eqref{on-diag}, respectively. Let~$z$ be the site where the maximum~$R_n$ is achieved and denote
\begin{equation}
%\label{}
\scrO_n:=\bigl\{x\in\scrC_{\infty,\alpha}\colon|x|\le n,\,|\psi_\omega(x)|\ge\tfrac12\epsilon n\bigr\}.
\end{equation}
Let~$Y=(Y_t)$ be a continuous-time random walk on~$\scrC_{\infty,\alpha}$ with expectation for the walk started at~$z$ denoted by~$E_{\omega,z}$. Define the stopping time
\begin{equation}
%\label{}
S_n:=\inf\bigl\{t\ge0\colon |Y_t-z|\ge2n\bigr\}
\end{equation}
and note that, in light of Proposition~\ref{prop-percolace}, we have $|Y_{t\wedge S_n}-z|\le 3n$ for all~$t>0$ provided $n\ge n_1(\omega)$ where~$n_1(\omega)<\infty$ a.s. The harmonicity of~$x\mapsto x+\psi_\omega(x)$ and the Optional Stopping Theorem yield
\begin{equation}
\label{M-exp}
R_n\le E_{\omega,z}\bigl|\psi_\omega(Y_{t\wedge S_n})+Y_{t\wedge S_n}-z\bigr|.
\end{equation}
Restricting to $t$ satisfying
\begin{equation}
\label{t-lower}
t\ge b_{4n},
\end{equation}
where~$b_n=o(n^2)$ is the sequence in part~(4) of Theorem~\ref{thm-sublinear}, we will now estimate the expectation separately on $\{S_n<t\}$ and $\{S_n\ge t\}$.

On the event~$\{S_n<t\}$, the absolute value in the expectation can simply be bounded by $R_{3n}+3n$. To estimate the probability of~$S_n<t$ we decompose according to whether $|Y_{2t}-z|\ge\frac32n$ or not. For the former, \eqref{t-lower} and \eqref{diffusive} imply
\begin{equation}
%\label{}
P_{\omega,z}\bigl(|Y_{2t}-z|\ge\tfrac32n\bigr)\le\frac{E_{\omega,z}|Y_{2t}-z|}{\frac32 n}\le
\frac 23C_1\frac{\sqrt{2t}}n.
\end{equation}
For the latter we invoke the inclusion
\begin{equation}
%\label{}
\bigl\{|Y_{2t}-z|\le\tfrac32n\bigr\}\cap\{S_n< t\}\subset\bigl\{|Y_{2t}-Y_{S_n}|\ge\tfrac12n\bigr\}\cap\{S_n< t\}
\end{equation}
and note that $2t-S_n\in[t,2t]$, \eqref{t-lower} and \eqref{diffusive} give us similarly
\begin{equation}
%\label{}
P_{\omega,x}\bigl(|Y_s-x|\ge \ffrac n2\bigr)\le \frac 2nC_1\sqrt{2t}\qquad\text{ when }x:=Y_{S_n}\text{ and }s:=2t-S_n.
\end{equation}
From the Strong Markov Property we thus conclude that this serves also as a bound for $P_{\omega,z}(S_n<t, |Y_{2t}-z|\ge\tfrac32n)$. Combining both parts and using $\frac83\sqrt2\le4$ we thus have
\begin{equation}
%\label{}
P_{\omega,z}(S_n< t)\le\frac{4C_1\sqrt{t}}n.
\end{equation}
The $S_n<t$ part of the expectation \eqref{M-exp} is bounded by $R_{3n}+3n$ times as much.

On the event $\{S_n\ge t\}$, the expectation in \eqref{M-exp} is bounded by
\begin{equation}
E_{\omega,z}\bigl(|\psi_\omega(Y_t)|\1_{\{S_n\ge t\}}\bigr)+ E_{\omega,z}|Y_t-z|.
\end{equation}
The second term on the right-hand side is then less than~$C_1\sqrt t$ provided~$t\ge b_n$. The first term is estimated depending on whether~$Y_t\not\in\scrO_{2n}$ or not:
\begin{equation}
%\label{}
E_{\omega,z}\bigl(|\psi_\omega(Y_t)|\1_{\{S_n\ge t\}}\bigr)
\le \frac12\epsilon n+R_{3n}P_{\omega,z}(Y_t\in\scrO_{2n}).
\end{equation}
For the probability of $Y_t\in\scrO_{2n}$ we get
\begin{equation}
%\label{}
P_{\omega,z}(Y_t\in\scrO_{2n})=\sum_{x\in\scrO_{2n}}P_{\omega,z}(Y_t=x)
\end{equation}
which, in light of the Cauchy-Schwarz estimate
\begin{equation}
%\label{}
P_{\omega,z}(Y_t=x)^2\le
P_{\omega,z}(Y_t=z)P_{\omega,x}(Y_t=x)
\end{equation}
and the definition of~$C_2$, is further estimated by
\begin{equation}
%\label{}
P_{\omega,z}(Y_t\in\scrO_{2n})\le C_2
\,\frac{|\scrO_{2n}|}{t^{d/2}}.
\end{equation}
From the above calculations we conclude that
\begin{equation}
%\label{}
R_n\le (R_{3n}+3n)\frac{4C_1\sqrt{t}}n+C_1\sqrt t+\frac12\epsilon n+R_{3n}C_2
\,\frac{|\scrO_{2n}|}{t^{d/2}}.
\end{equation}
Since $|\scrO_{2n}|=o(n^d)$ as~$n\to\infty$, by \eqref{on-average} we can choose~$t:=\xi n^2$ with~$\xi>0$ so small that \eqref{t-lower} applies and \eqref{mn-bd} holds for the given~$\epsilon$ and~$\delta$ once~$n$ is sufficiently large.
\end{proofsect}

We now proceed to prove convergence of the random walk~$X=(X_n)$ to Brownian motion. Most of the ideas are drawn directly from Berger and Biskup~\cite{BB} so we stay rather brief. We will frequently work on the truncated infinite component~$\scrC_{\infty,\alpha}$ and the corresponding restriction of the random walk; cf~\twoeqref{T-formula}{generator}. We assume throughout that~$\alpha$ is such that \twoeqref{infinite}{finite}~hold. 

\begin{lemma}
\label{lemma-harmonic}
Let~$\chi$ be the corrector on~$\scrC_\infty$. Then $\varphi_\omega(x):=x+\chi(\omega,x)$ is harmonic for the random walk observed only on~$\scrC_{\infty,\alpha}$, i.e.,
\begin{equation}
\LL_\omega^{(\alpha)}\varphi_\omega(x)=0,\qquad\forall x\in\scrC_{\infty,\alpha}.
\end{equation}
\end{lemma}

\begin{proofsect}{Proof}
We have
\begin{equation}
(\LL_\omega^{(\alpha)}\varphi_\omega)(x)=E_{\omega,x}\bigl(\varphi_\omega(X_{T_1})\bigr)-\varphi_\omega(x)
\end{equation}
But~$X_n$ is confined to a finite component of~$\scrC_\infty\setminus\scrC_{\infty,\alpha}$ for~$n\in[0,T_1]$, and so~$\varphi_\omega(X_n)$ is bounded. Since~$(\varphi_\omega(X_n))$ is a martingale and~$T_1$ is an a.s.\ finite stopping time, the Optional Stopping Theorem tells us $E_{\omega,x}\varphi_\omega(X_{T_1})=\varphi_\omega(x)$.
\end{proofsect}

Next we recall the proof of sublinearity of the corrector along coordinate directions:

\begin{lemma}
\label{lemma-along-axes}
For~$\omega\in\{0\in\scrC_{\infty,\alpha}\}$, let~$(x_n(\omega))_{n\in\Z}$ mark the intersections of~$\scrC_{\infty,\alpha}$ and one of the coordinate axis so that $x_0(\omega)=0$. Then
\begin{equation}
\lim_{n\to\infty}\frac{\chi(\omega,x_n(\omega))}n=0,\qquad\BbbP_\alpha\text{\rm-a.s.}
\end{equation}
\end{lemma}

\begin{proofsect}{Proof}
Let~$\tau_x$ be the ``shift by~$x$'' on~$\Omega$ and let~$\sigma(\omega):=\tau_{x_1(\omega)}(\omega)$ denote the ``induced'' shift. Standard arguments (cf~\cite[Theorem~3.2]{BB}) prove that $\sigma$ is~$\BbbP_\alpha$ preserving and ergodic. Moreover, 
\begin{equation}
\E_\alpha\bigl(\textd_\omega^{(\alpha)}(0,x_1(\omega))^p\bigr)<\infty,\qquad p<\infty,
\end{equation}
by \cite[Lemma~4.3]{BB} (based on Antal and Pisztora~\cite{Antal-Pisztora}). Define $\Psi(\omega)\,:=\chi(\omega,x_1(\omega))$. Theorem~\ref{thm3.1}(5) tells us that
\begin{equation}
\Psi\in L^1(\BbbP_\alpha)\quad\text{and}\quad \E_\alpha\Psi(\omega)=0.
\end{equation}
But the gradient property of~$\chi$ implies
\begin{equation}
\frac{\chi(\omega,x_n(\omega))}n=\frac1n\sum_{k=0}^{n-1}\Psi\circ\sigma^k(\omega)
\end{equation}
and so the left-hand side tends to zero a.s.\ by the Pointwise Ergodic Theorem.
\end{proofsect}

We will also need sublinearity of the corrector on average:

\begin{lemma}
\label{lemma-on-average}
For each~$\epsilon>0$ and~$\BbbP_\alpha$-a.e.~$\omega$:
\begin{equation}
%\label{}
\lim_{n\to\infty}\frac1{n^d}\sum_{\begin{subarray}{c}
x\in\scrC_{\infty,\alpha}\\|x|\le n
\end{subarray}}
\1_{\{|\chi(\omega,x)|\ge\epsilon n\}}=0.
\end{equation}
\end{lemma}

\begin{proofsect}{Proof}
This follows from Lemma~\ref{lemma-along-axes} exactly as~\cite[Theorem~5.4]{BB}.
\end{proofsect}

\begin{remark}
The proof of \cite[Theorem~5.4]{BB} makes a convenient use of separate ergodicity (i.e., that with respect to shifts only in one of the coordinate directions). This is sufficient for i.i.d.\ environments as considered in the present situation. However, it is not hard to come up with a modification of the proof that covers general ergodic environments as well (Biskup and Deuschel~\cite{Biskup-Deuschel}). 
\end{remark}

Finally, we will assert the validity of the bounds on the return probability and expected displacement of the walk from Theorem~\ref{thm-sublinear}:

\begin{lemma}
\label{lemma-diffusive}
Let~$(Y_t)$ denote the continuous-time random walk on~$\scrC_{\infty,\alpha}$. Then the diffusive bounds \twoeqref{diffusive}{on-diag} hold for~$\BbbP_\alpha$-a.e.~$\omega$.
\end{lemma}

We will prove this lemma at the very end of Sect.~\ref{sec5}.

\begin{proofsect}{Proof of Theorem~\ref{thm-main}}
Let~$\alpha$ be such that \twoeqref{infinite}{finite} hold and let~$\chi$ denote the corrector on~$\scrC_\infty$ as constructed in Theorem~\ref{thm3.1}. The crux of the proof is to show that $\chi$ grows sublinearly with~$x$, i.e., $\chi(\omega,x)=o(|x|)$ a.s.
 
By Lemmas~\ref{lemma-harmonic} and \ref{lemma-on-average}, Theorem~\ref{thm3.1}(4) and Lemma~\ref{lemma-diffusive}, the corrector satisfies the conditions of Theorem~\ref{thm-sublinear}. It follows that $\chi$ is sublinear on~$\scrC_{\infty,\alpha}$ as stated in~\eqref{sublinear}. However, by \eqref{exponential} the largest component of $\scrC_\infty\setminus\scrC_{\infty,\alpha}$ in a box~$[-2n,2n]$ is less than~$C\log n$ in diameter, for some random but finite~$C=C(\omega)$. Invoking the harmonicity of~$\varphi_\omega$ on~$\scrC_\infty$, the Optional Stopping Theorem~gives
\begin{equation}
\label{corrector-sublinear}
\max_{\begin{subarray}{c}
x\in\scrC_\infty\\|x|\le n
\end{subarray}}
\bigl|\chi(\omega,x)\bigr|
\le \max_{\begin{subarray}{c}
x\in\scrC_{\infty,\alpha}\\|x|\le n
\end{subarray}}
\bigl|\chi(\omega,x)\bigr|+2C(\omega)\log(2n),
\end{equation}
whereby we deduce that $\chi$ is sublinear on~$\scrC_\infty$ as well.

Having proved the sublinearity of~$\chi$ on~$\scrC_\infty$, we proceed as in the~$d=2$ proof of~\cite{BB}. Let $\varphi_\omega(x):=x+\chi(\omega,x)$ and abbreviate $M_n:=\varphi_\omega(X_n)$. Fix~$\hatv\in\R^d$ and define
\begin{equation}
f_K(\omega):=E_{\omega,0}\bigl((\hatv\cdot M_1)^2\1_{\{|\hatv\cdot M_1|\ge K\}}\bigr).
\end{equation}
By Theorem~\ref{thm3.1}(3), $f_K\in L^1(\Omega,\scrF,\BbbP_0)$ for all~$K$. Since the Markov chain on environments, $n\mapsto\tau_{X_n}(\omega)$, is ergodic (cf~\cite[Section~3]{BB}), we thus have
\begin{equation}
\label{5.28-eq}
\frac1n\sum_{k=0}^{n-1}f_K\circ\tau_{X_k}(\omega)\,\underset{n\to\infty}\longrightarrow\,\E_0 f_K,
\end{equation}
for $\BbbP_0$-a.e.~$\omega$ and $P_{\omega,0}$-a.e.\ path~$X=(X_k)$ of the random walk. Using this for~$K:=0$ and~$K:=\epsilon\sqrt n$ along with the monotonicity of~$K\mapsto f_K$ verifies the conditions of the Lindeberg-Feller Martingale Functional CLT (e.g., Durrett~\cite[Theorem~7.7.3]{Durrett}). Thereby we conclude that the random continuous function
\begin{equation}
\label{linear-M}
t\mapsto \frac1{\sqrt n}\bigl(\hatv\cdot M_{\lfloor nt\rfloor}+(nt-\lfloor nt\rfloor)\,\hatv\cdot(M_{\lfloor nt\rfloor+1}-M_{\lfloor nt\rfloor})\bigr)
\end{equation}
converges weakly to Brownian motion with mean zero and covariance
\begin{equation}
\E_0 f_0=\E_0E_{\omega,0}\bigl((\hatv\cdot M_1)^2\bigr).
\end{equation}
This can be written as~$\hatv\cdot D\hatv$ where~$D$ is the matrix with coefficients
\begin{equation}
D_{i,j}:=\E_0 E_{\omega,0}\bigl((\hate_i\cdot M_1)(\hate_j\cdot M_1)\bigr).
\end{equation}
Invoking the Cram\'er-Wold device (e.g., Durrett~\cite[Theorem~2.9.5]{Durrett}) and the fact that continuity of a stochastic process in~$\R^d$ is implied by the continuity of its~$d$ one-dim\-en\-sional projections we get that the linear interpolation of~$t\mapsto M_{\lfloor nt\rfloor}/\sqrt n$ scales to $d$-dimensional Brownian motion with covariance matrix~$D$. The sublinearity of the corrector then ensures, as in~\cite[(6.11--6.13)]{BB}, that
\begin{equation}
X_n-M_n=\chi(\omega,X_n)=o(|X_n|)=o(|M_n|)=o(\sqrt n),
\end{equation}
and so the same conclusion applies to~$t\mapsto B_n(t)$ in \eqref{1.5a}.

The reflection symmetry of~$\BbbP_0$ forces~$D$ to be diagonal; the rotation symmetry then ensures that~$D=(\ffrac1d)\sigma^2\1$ where
\begin{equation}
%\label{}
\sigma^2:=\E_0E_{\omega,0}|M_1|^2
\end{equation}
To see that the limiting process is not degenerate to zero we note that if we had $\sigma=0$ then $\chi(\cdot,x)=-x$ would hold a.s.\ for all~$x\in\Z^d$. But that is impossible since, as we proved above,~$x\mapsto\chi(\cdot,x)$ is sublinear~a.s.
\end{proofsect}

\begin{remark}
Note that, unlike the proofs in~\cite{Sidoravicius-Sznitman,BB,Mathieu-Piatnitski}, the above line of argument does not require a separate proof of tightness. In our approach, this comes rather automatically for the deformed random walk $\varphi_\omega(X_n)$ --- via the (soft) stationarity argument~\eqref{5.28-eq} and the Martingale Functional CLT. Sublinearity of the corrector then extends it readily to the original random walk.
\end{remark}

\begin{remark}
We also wish to use the opportunity to correct an erroneous argument from~\cite{BB}. There, at the end of the proof of Theorem~6.2 it is claimed that the expectation $\E_0 E_{\omega,0}(X_1\cdot\chi(X_1,\omega))$ is zero. Unfortunately, this is false. In fact, we have
\begin{equation}
%\label{}
\E_0 E_{\omega,0}\bigl(X_1\cdot\chi(X_1,\omega)\bigr)=-\E_0 E_{\omega,0}\bigl|\chi(X_1,\omega)\bigr|^2<0.
\end{equation}
where the strict inequality assumes that~$\BbbP$ is non-degenerate. This shows 
\begin{equation}
%\label{}
\E_0 E_{\omega,0}|M_1|^2=\E_0 E_{\omega,0}|X_1|^2-\E_0 E_{\omega,0}\bigl|\chi(X_1,\omega)\bigr|^2<\E_0 E_{\omega,0}|X_1|^2.
\end{equation}
Thus, once~$\BbbP$ is non-degenerate, the diffusion constant of the limiting Brownian motion is strictly \emph{smaller} than the variance of the first step.

A consequence of the above error for the proof of Theorem~6.2 in~\cite{BB} is that it invalidates one of the three listed arguments to prove that the limiting Brownian motion is non-degenerate. Fortunately, the remaining two arguments are correct.
\end{remark}

\section{Heat kernel and expected distance}
\label{sec5}
\noindent
Here we will derive the bounds \twoeqref{diffusive}{on-diag} and thus establish Lemma~\ref{lemma-diffusive}. Most of the derivation will be done for a general countable-state Markov chain; we will specialize to random walk among i.i.d.\ conductances at the very end of this section. The general ideas underlying these derivations are fairly standard and exist, in some form, in the literature. A novel aspect is the way we control the non-uniformity of volume-growth caused by local irregularities of the underlying graph; cf~\eqref{C-vol} and Lemma~\ref{lemma-QM}(1). A well informed reader may nevertheless wish to read only the statements of Propositions~\ref{lemma-upper} and~\ref{thm-distance} and then pass directly to the proof of Lemma~\ref{lemma-diffusive}.

\smallskip
Let~$V$ be a countable set and let~$(a_{xy})_{x,y\in V}$ denote the collection of positive numbers with the following properties: For all~$x,y\in V$,
\begin{equation}
%\label{}
a_{xy}=a_{yx}\quad\text{and}\quad\pi(x)\,:=\sum_{y\in V}a_{xy}<\infty.
\end{equation}
Consider a continuous time Markov chain $(Y_t)$ on~$V$ with the generator
\begin{equation}
%\label{}
(\LL f)(x):=\frac1{\pi(x)}\sum_{y\in V}a_{xy}\bigl[f(y)-f(x)\bigr].
\end{equation}
We use~$P^x$ to denote the law of the chain started from~$x$, and~$E^x$ to denote the corresponding expectation.
Consider a graph~$G=(V,E)$ where~$E$ is the set of all pairs~$(x,y)$ such that~$a_{xy}>0$. Let~$\textd(x,y)$ denote the distance between~$x$ and~$y$ as measured on~$G$.

For each~$x\in V$, let~$B_n(x):=\{y\in V\colon\textd(x,y)\le n\}$. If~$\Lambda\subset V$, we use $Q(\Lambda,\Lambda^\cc)$ to denote the sum
\begin{equation}
%\label{}
Q(\Lambda,\Lambda^\cc):=\sum_{x\in\Lambda}\sum_{y\in\Lambda^\cc}a_{xy}.
\end{equation}
Suppose that there are constants~$d>0$ and $\nu\in(0,\ffrac12)$ such that, for some~$a>0$,
\begin{equation}
\label{C-vol}
\Cvol (x,a)\,:=\sup_{0<s\le a}\,\Bigl[\,s^d\sum_{y\in V}\pi(y)\texte^{-s\textd(x,y)}\Bigr]<\infty
\end{equation}
and
\begin{equation}
\label{C-iso}
\Ciso (x)\,:=\inf_{n\ge1}\inf\biggl\{\,\frac{Q(\Lambda,\Lambda^\cc)}{\pi(\Lambda)^{\frac{d-1}d}}\colon
\Lambda\subset B_{2n}(x),\,\pi(\Lambda)\ge n^\nu,\,\Lambda\text{ connected}\biggr\}>0.
\end{equation}
Let~$V(\epsilon)\subset V$ denote the set of all~$x\in V$ that are connected to infinity by a self-avoiding path $(x_0=x,x_1,\dots)$ with $a_{x_i,x_{i+1}}\ge\epsilon$ for all~$i\ge0$. Suppose that
\begin{equation}
\label{c-mini}
a_\star\,:=\sup\bigl\{\epsilon>0\colon V(\epsilon)=V\bigr\}>0.
\end{equation}
(Note that this does not require~$a_{xy}$ be bounded away from zero.)

\smallskip
The first observation is that the heat-kernel, defined by
\begin{equation}
%\label{}
q_t(x,y)\,:=\frac{P^x(Y_t=y)}{\pi(y)},
\end{equation}
can be bounded in terms of the isoperimetry constant~$\Ciso (x)$. Bounds of this form are well known and have been derived by, e.g., Coulhon, Grigor'yan and Pittet~\cite{CGP} for heat-kernel on manifolds, and by Lov\'asz and Kannan~\cite{Lovasz-Kannan}, Morris and Peres~\cite{Morris-Peres} and Goel, Montenegro and Tetali~\cite{GM-Tetali} in the context of countable-state Markov chains. We will use the formulation for infinite graphs developed in Morris and Peres~\cite{Morris-Peres}.

\begin{proposition}
\label{lemma-upper}
There exists a constant~$c_1\in(1,\infty)$ depending only on~$d$ and~$a_\star$ such that for $t(x)\,:=c_1[\log (\Ciso (x)\vee c_1)]^{\frac1{1-2\nu}}$ we have
\begin{equation}
\label{q-uniform}
\sup_{z\in B_t(x)}\,\sup_{y\in V}\,q_t(z,y)\le c_1\frac{\Ciso (x)^{-d}}{t^{d/2}},
\qquad t\ge t(x).
\end{equation}
\end{proposition}

\begin{proofsect}{Proof}
We will first derive the corresponding bound for the discrete-time version of~$(Y_t)$.
Let $\cmss P(x,y):=a_{xy}/\pi(x)$ and define $\hat{\cmss P}:=\frac12(1+\cmss P)$. Let
\begin{equation}
%\label{}
\hat q_n(x,y):=\frac{\hat{\cmss P}^n(x,y)}{\pi(y)}
\end{equation}
We claim that, for some absolute constant~$c_1$ and any~$z\in B_n(x)$,
\begin{equation}
\label{discrete-heat}
\hat q_n(z,y)\le c_1\frac{\Ciso(x)^{-d}}{n^{d/2}},\qquad n\ge t(x).
\end{equation}
To this end, let $\widehat Q$ be the object~$Q$ for the Markov chain~$\hat{\cmss P}$ and let
\begin{equation}
\label{phi-def}
\phi(r):=\inf\Bigl\{\frac{\widehat Q(\Lambda,\Lambda^\cc)}{\pi(\Lambda)}\colon\pi(\Lambda)\le r,\,\Lambda\subset B_{2n}(x)\Bigr\}.
\end{equation}
We claim that Theorem~2 of Morris and Peres~\cite{Morris-Peres} then implies that, for any $\epsilon$ that satisfies $0<\epsilon<[\pi(x)\wedge\pi(y)]^{-1}$ and
\begin{equation}
%\label{}
n\ge1+\int_{4(\pi(z)\wedge\pi(y))}^{4/\epsilon} \frac{4\textd  r}{r\phi(r)^2},
\end{equation}
we have $\hat q_n(z,y)\le\epsilon$. 
To see this we have to check that the restriction~$\Lambda\subset B_{2n}(x)$ in the definition of~$\phi(r)$, which is absent from the corresponding object in~\cite{Morris-Peres}, causes no harm. First note that the Markov chain started at $z\in B_n(x)$ will not leave $B_{2n}(x)$ by time~$n$. Thus, we can modify the chain outside~$B_{2n}(x)$ arbitrarily. It is easy to come up with a modification that will effectively reduce the infimum in \eqref{phi-def} to sets inside~$B_{2n}(x)$.

It is well known (and easy to check) that the infimum in \eqref{phi-def} can be restricted to connected sets~$\Lambda$. Then \twoeqref{C-iso}{c-mini} give us
\begin{equation}
\label{phin}
\phi(r)\ge\frac12\bigl(\Ciso (x)r^{-1/d}\wedge a_\star n^{-\nu}\bigr)
\end{equation}
where the extra half arises due the consideration of time-delayed chain $\hat{\cmss P}$.
The two regimes cross over at $r_n\,:=(\Ciso (x)/a_\star)^dn^{d\nu}$; the integral is thus bounded by
\begin{equation}
\label{integral-comp}
\int_{4(\pi(z)\wedge\pi(y))}^{4/\epsilon} \frac{4\textd  r}{r\phi(r)^2}
\le4\frac{n^{2\nu}}{a_\star^2}\log\Bigl(\frac{r_n}{4a_\star}\Bigr)+2d\Ciso (x)^{-2}\Bigl(\frac4\epsilon\Bigr)^{2/d}.
\end{equation}
The first term splits into a harmless factor of order $n^{2\nu}\log n=o(n)$ and a term proportional to $n^{2\nu}\log \Ciso (x)$. This is~$O(n)$ by~$n\ge t(x)$ where the (implicit) constant can be made as small as desired by choosing~$c_1$ sufficiently large. Setting 
\begin{equation}
%\label{}
\epsilon:=c[\Ciso(x)^2n]^{-d/2}
\end{equation}
we can thus adjust the constant~$c$ in such a way that \eqref{integral-comp} is less than~$n-1$ for all $n\ge t(x)$. Thereby \eqref{discrete-heat} follows.

To extend the bound \eqref{discrete-heat} to continuous time, we note that $\LL=2(\hat{\cmss P}-1)$. Thus if~$N_t$ is Poisson with parameter~$2t$, then
\begin{equation}
\label{discrete-continuous}
q_t(z,y)=E\hat q_{N_t}(z,y).
\end{equation}
But~$P(N_t\le \tfrac32 t\,\text{or}\,N_t\ge3t)$ is exponentially small in~$t$ and, in particular, much smaller than \eqref{q-uniform} for $t\ge c_1\log \Ciso(x)$ with~$c_1$ sufficiently large. As~$q_t\le(a_\star)^{-1}$, the~$N_t\not\in(\tfrac32t,3t)$ portion of the expectation in \eqref{discrete-continuous} is thus negligible. Once~$N_t$ is constrained to the interval $(t,3t)$ the uniform bound \eqref{discrete-heat} implies \eqref{q-uniform}.
\end{proofsect}

Our next item of business is a diffusive bound on the expected (graph-theoretical) distance traveled by the walk~$Y_t$ by time~$t$. As was noted by Bass~\cite{Bass} and Nash~\cite{Nash}, this can be derived from the above uniform bound on the heat-kernel assuming regularity of the volume growth. Our proof is an adaptation of an argument of Barlow~\cite{Barlow}.

\begin{proposition}
\label{thm-distance}
There exist constants~$c_2=c_2(d)$ and $c_3=c_3(d)$  such that the following holds: Let~$x\in V$ and suppose $A>0$ and $t(x)>1$ are numbers for which
\begin{equation}
\label{5.15}
\sup_{y\in V}q_t(x,y)\le\frac{A}{t^{d/2}},\qquad t\ge t(x),
\end{equation}
holds and let $T(x):=\frac1d(Aa_\star)^{-4/d}\vee [t(x)\log t(x)]$. Then
\begin{equation}
E^x\textd(x,Y_t) \le A'(x,t)\sqrt t,\qquad t\ge T(x),
\end{equation}
with~$A'(x,t):=c_2+c_3[\log A+\Cvol(x,t^{-1/2})]$.
\end{proposition}

Much of the proof boils down to the derivation of rather inconspicuous but deep relations (discovered by Nash~\cite{Nash}) between the following quantities:
\begin{equation}
%\label{}
M(x,t)\,:=E^x\textd(x,Y_t) =\sum_y \pi(y)q_t(x,y)\textd(x,y)
\end{equation}
and
\begin{equation}
%\label{}
Q(x,t)\,:=-E^x \log q_t(x,Y_t) = -\sum_y \pi(y)q_t(x,y)\log q_t(x,y).
\end{equation}
Note that~$q_t(x,\cdot)\le(a_\star)^{-1}$ implies $Q(x,t)\ge\log a_\star$.

\begin{lemma}
\label{lemma-QM}
For all~$t\ge0$ and all~$x\in V$,
\settowidth{\leftmargini}{(111a)}
\begin{enumerate}
\item[(1)]
$M(x,t)^d\ge \exp\{-1-\Cvol (x,M(x,t)^{-1})+Q(x,t)\}$,
\item[(2)]
$M'(x,t)^2\le Q'(x,t)$.
\end{enumerate}
\end{lemma}

\begin{proof}
Part~(2) is identical to the proof of Lemma~3.3 in Barlow~\cite{Barlow} so we only need to prove part~(1). We will follow the proof of \cite[Lemma~3.3]{Barlow} except that we need to keep track of $\Cvol (x)$ in the calculation. Pick two numbers~$a>0$ and~$b\in\R$ and note that
\begin{equation}
%\label{}
u\log u+\lambda u\ge-\texte^{-\lambda-1}.
\end{equation}
This implies
\begin{equation}
%\label{}
-Q(x,t)+aM(x,t)+b\ge-\sum_y\pi(y)\texte^{-b-1-a\textd(x,y)}.
\end{equation}
Using the definition of~$\Cvol(x,a)$ and bounding~$\texte^{-1}\le1$ we get
\begin{equation}
%\label{}
-Q(x,t)+aM(x,t)+b\ge -\Cvol (x,a)\,\texte^{-b}a^{-d}.
\end{equation}
Now set~$\texte^{-b}:=a^d$ with~$a:=M(x,t)^{-1}$ to get the result.
\end{proof}

These bounds imply the desired diffusive estimate on~$M(x,t)$:

\begin{proofsect}{Proof of Proposition~\ref{thm-distance}}
Suppose without loss of generality that~$M(x,t)\ge\sqrt t$, because otherwise there is nothing to prove.
We follow the proof of \cite[Proposition~3.4]{Barlow}. The key input is provided by the inequalities in Lemma~\ref{lemma-QM}. Define the function
\begin{equation}
L(t):=\frac1d\Bigl(Q(x,t)+\log A-\frac d2\log t\Bigr)
\end{equation}
and note that~$L(t)\ge0$ for~$t\ge t(x)$. Let
\begin{equation}
%\label{}
t_0:=(Aa_\star)^{-2/d}\vee\sup\bigl\{t\ge0\colon L(t)\le0\bigr\}.
\end{equation}
We claim that~$M(x,t_0)\le\sqrt{dT(x)}$. Indeed, when $t_0=(Aa_\star)^{-2/d}$ then this follows by
\begin{equation}
M(x,t_0)\le t_0=(Aa_\star)^{-2/d}\le\sqrt{dT(x)}
\end{equation}
due to our choice of~$T(x)$. On the other hand, when~$t_0>(Aa_\star)^{-2/d}$ we use Lemma~\ref{lemma-QM}(2), the Fundamental Theorem of Calculus and the Cauchy-Schwarz inequality to derive
\begin{equation}
%\label{}
M(x,t_0)\le\sqrt{t_0}\bigl[Q(x,t_0)-Q(x,0)\bigr]^{1/2}.
\end{equation}
Since $Q(x,0)\ge\log a_\star$ and $L(t_0)=0$ by continuity, we have
\begin{equation}
M(x,t_0)\le\sqrt{t_0}\,\Bigl(\,\frac d2\log t_0-\log A-\log a_\star\Bigr)^{1/2}
\le\sqrt{dt_0\log t_0}
\end{equation}
where we used that~$t_0\ge(Aa_\star)^{-2/d}$ implies $\log A+\log a_\star\ge-\frac d2\log t_0$. Since this implies~$t_0\ge1$, the condition $t_0\le t(x)$ shows that the right-hand side is again less than~$\sqrt{dT(x)}$.

For~$t\ge t_0$ we have~$L(t)\ge0$. Lemma~\ref{lemma-QM}(2) yields
\begin{equation}
\label{MMbd}
\begin{aligned}
M(x,t)-M(x,t_0)&\le\sqrt d\,\int_{t_0}^t \Bigl(\frac1{2s}+L'(s)\Bigr)^{1/2}\textd s
\\
&\le\sqrt d\,\int_{t_0}^t \Bigl(\frac1{\sqrt{2s}}+L'(s)\sqrt{\ffrac s2}\,\Bigr)\textd s
\le\sqrt{2dt}+L(t)\sqrt{dt},
\end{aligned}
\end{equation}
where we used integration by parts and the positivity of~$L$ to derive the last inequality. Now put this together with~$M(x,t_0)\le\sqrt{dt}$ and apply Lemma~\ref{lemma-QM}(1), noting that $\Cvol(z,M(x,t)^{-1})\le\Cvol(z,t^{-1/2})$ by the assumption $M(x,t)\ge\sqrt t$. Dividing out an overall factor~$\sqrt t$, we thus get
\begin{equation}
\bigl[A\texte^{\Cvol(x,t^{-1/2})}\bigr]^{-1/d}\texte^{-1/d+L(t)}\le3\sqrt{d}+\sqrt{d}\,L(t).
\end{equation}
This implies
\begin{equation}
%\label{}
L(t)\le \tilde c_2+\tilde c_3\bigl[\log A+\Cvol (x,t^{-1/2})\bigr]
\end{equation}
for some constants~$\tilde c_2$ and $\tilde c_3$ depending only on~$d$. Plugging this in \eqref{MMbd}, we get the desired claim.
\end{proofsect}

We are now finally ready to complete the proof of our main theorem:

\begin{proofsect}{Proof of Lemma~\ref{lemma-diffusive}}
We will apply the above estimates to obtain the proof of the bounds \twoeqref{diffusive}{on-diag}. We use the following specific choices
\begin{equation}
%\label{}
V:=\scrC_{\infty,\alpha},\quad a_{xy}:=\hat\omega_{xy},\quad\pi(x):=2d,\quad\text{and} \quad b_n:=n.
\end{equation}
As $a_\star\ge\alpha$, all required assumptions are satisfied.  

To prove \eqref{on-diag}, we note that by Lemma~3.4 of Berger, Biskup, Hoffman and Kozma~\cite{BBHK} (based on the isoperimetric inequality for the supercritical bond-percolation cluster, cf\ Benjamini and Mossel~\cite{Benjamini-Mossel}, Rau~\cite[Proposition~1.2]{Rau} or \cite[Section~5]{BBHK}) we have $\Ciso (0)>0$ a.s. Hence, Proposition~\ref{lemma-upper} ensures that, for all~$z\in\scrC_{\infty,\alpha}$ with $|z|\le t$,
\begin{equation}
%\label{}
t^{d/2}P_{\omega,z}(Y_t=z)\le2dc_1\Ciso (0)^{-d}
\end{equation}
provided~$t$ exceeds some~$t_1$ depending on~$\Ciso (0)$. From here \eqref{on-diag} immediately follows. 

To prove \eqref{diffusive}, we have to show that, a.s.,
\begin{equation}
\label{Cvol-bd}
\sup_{n\ge1}\,\,\max_{\begin{subarray}{c}
z\in\scrC_{\infty,\alpha}\\|z|\le n
\end{subarray}}\,\,\sup_{t\ge n}\,\,
\Cvol(z,t^{-1/2})<\infty.
\end{equation}
To this end we note that Lemma~\ref{lemma-compare} implies that there is a.s.\ finite $C=C(\omega)$ such that for all~$z,y\in\scrC_{\infty,\alpha}$ with~$|z|\le n$ and~$|z-y|\ge C\log n$,
\begin{equation}
\label{dist-comp}
\textd(z,y)\ge\varrho|z-y|.
\end{equation}
It follows that, once~$\ffrac1a>C\log n$, for every~$z\in\scrC_{\infty,\alpha}$ with~$|z|\le n$ we have
\begin{equation}
%\label{}
\sum_{y\in\scrC_{\infty,\alpha}}
\texte^{-a\textd(z,y)}
\le c_4 a^{-d}+\sum_{\begin{subarray}{c}
y\in\scrC_{\infty,\alpha}\\|y-z|\ge 1/a
\end{subarray}}
\texte^{-a\varrho|z-y|}\le c_5 a^{-d},
\end{equation}
where~$c_4$ and~$c_5$ are constants depending on~$d$ and~$\varrho$. Since~$\ffrac1a=t^{1/2}\ge\sqrt n\gg\log n$, \eqref{Cvol-bd} follows.

Once we have the uniform bound \eqref{Cvol-bd}, as well as the uniform bound \eqref{5.15}, Proposition~\ref{thm-distance} yields the a.s.\ inequality
\begin{equation}
%\label{}
\sup_{n\ge1}\,\,\max_{\begin{subarray}{c}
z\in\scrC_{\infty,\alpha}\\|z|\le n
\end{subarray}}\,\,\sup_{t\ge n}\,\,
\frac{E_{\omega,z}\textd(z,Y_t)}{\sqrt t}<\infty.
\end{equation}
To convert $\textd(z,Y_t)$ into~$|z-Y_t|$ in the expectation, we invoke \eqref{dist-comp} one more time.
\end{proofsect}

\section*{Acknowledgments}
\noindent 
The research of M.B.~was supported by the NSF grant~DMS-0505356. The authors wish to thank N.~Berger, M.~Barlow, J.-D.~Deuschel, C.~Hoffman, G.~Kozma, Y.~Peres and G.~Pete for discussions on this problem.

%%%%%%%%%%%%%%%%%%%%%%%%%%%%%%%%%%%%%%%%%%%%%%%%%%%%%
\end{document}